\documentclass[12pt]{article}
\title{A Nullstellensatz for amoebas}
\author{Kevin Purbhoo \\ University of British Columbia}

\usepackage[mathscr]{eucal}
\usepackage{latexsym}
\usepackage{amssymb}
\usepackage{amsthm}
\usepackage{amsmath}
\usepackage{epsfig}

\newcommand{\CC}{\mathbb{C}}
\newcommand{\ZZ}{\mathbb{Z}}
\newcommand{\RR}{\mathbb{R}}
\newcommand{\QQ}{\mathbb{Q}}
\newcommand{\PP}{\mathbb{P}}
\newcommand{\TT}{\mathbb{T}}
\newcommand{\Rtrop}{\mathbb{R}_\mathrm{trop}}
\newcommand{\Log}{\mathrm{Log}}
\newcommand{\Vol}{\mathrm{Vol}}
\newcommand{\Res}{\mathrm{Res}}
\newcommand{\Am}{\mathcal{A}}
\newcommand{\Spine}{\mathcal{S}}
\newcommand{\LopSpine}{\mathcal{LS}}
\newcommand{\LopAm}{\mathcal{LA}}
\newcommand{\SlopAm}{\mathcal{SA}}
\newcommand{\TropVar}{\mathcal{T}}
\newcommand{\Amcomp}{\RR^r \setminus \mathcal{A}}
\newcommand{\LopAmcomp}{\RR^r \setminus \mathcal{LA}}
\newcommand{\SlopAmcomp}{\RR^r \setminus \mathcal{SA}}
\newcommand{\boldz}{{\bf z}}
\newcommand{\boldw}{{\bf w}}
\newcommand{\boldx}{{\bf x}}
\newcommand{\boldy}{{\bf y}}
\newcommand{\bolda}{{\bf a}}
\DeclareMathAlphabet{\mathbsl}{OT1}{cmr}{bx}{sl}
\newcommand{\alta}{\mathbsl{a}}
\newcommand{\altz}{\mathbsl{z}}

\newcommand{\val}{{\bf val}}
\newcommand{\boldzeta}{\mbox{\boldmath{$\zeta$}}}
\newcommand{\smallboldzeta}{\mbox{\scriptsize\boldmath{$\zeta$}}}
\newcommand{\bdotsb}{\raisebox{1.5pt}{\text{\tiny$\{$}} \cdots \raisebox{1.5pt}{\text{\tiny$\}$}}}

\newcommand{\basering}{\CC[z_1,z_1^{-1}, \ldots, z_r,z_r^{-1}]}
\newcommand{\bigbasering}{\CC[z_1,z_1^{-1}, \ldots, z_r,z_r^{-1}, w_1, w_1^{-1}, \ldots, w_{r'}, w_{r'}^{-1}]}
\newcommand{\Kbasering}{K[z_1,z_1^{-1}, \ldots, z_r,z_r^{-1}]}
\newcommand{\gencyclres}[1]{\tilde{#1}_{n_1, \ldots, n_r}}
\newcommand{\cyclres}[1]{\tilde{#1}_n}
\newcommand{\widecyclres}[1]{\widetilde{#1}_n}
\newcommand{\projcyclres}[2]{\tilde{#1}_n^{#2,\smallboldzeta}}
\newcommand{\projf}[1]{f^{#1,\smallboldzeta}}
\newcommand{\projfprime}[1]{f^{#1,\smallboldzeta'}}
\newcommand{\maxdeg}{\mathrm{maxdeg}}
\newcommand{\mindeg}{\mathrm{mindeg}}
\newcommand{\mdeg}[1]{c_{#1} n^{D_{#1}}}
\newcommand{\vecj}{{\overrightarrow{j}}}
\newcommand{\veck}{{\overrightarrow{k}}}
\newcommand{\vecl}{{\overrightarrow{l}}}
\newcommand{\Proj}{\mathrm{Proj}}
\newcommand{\Spec}{\mathrm{Spec}}
\newcommand{\ind}{\mathrm{ind}}
\newcommand{\wt}{\mathrm{wt}}

\newtheorem{theorem}{Theorem}
\newtheorem*{theorem*}{Theorem}
\newtheorem{lemma}{Lemma}[section]
\newtheorem{corollary}[lemma]{Corollary}
\newtheorem{calculation}[lemma]{Calculation}
\newtheorem{proposition}[lemma]{Proposition}
\newtheorem{example}{Example}[section]
\newtheorem*{remark*}{Remark}
\newtheorem{remark}[example]{Remark}
\newtheorem*{definition*}{Definition}
\newtheorem{definition}[example]{Definition}

\newcommand{\xothertheoremname}{Theorem}
\newtheorem*{xothertheorem}{\xothertheoremname}
\newenvironment{othertheorem}[1]
{\renewcommand{\xothertheoremname}{#1}\begin{xothertheorem}}
{\end{xothertheorem}}

\numberwithin{equation}{section}

\begin{document}
\maketitle
\allowdisplaybreaks

\begin{abstract}
The amoeba of an affine algebraic variety $V \subset (\CC^*)^r$ is the image
of $V$ under the map $(z_1, \ldots, z_r) \mapsto
(\log|z_1|, \ldots, \log|z_r|)$.
We give a characterisation of
the amoeba based on the triangle inequality, which we call
`testing for {\em lopsidedness}'.  
We show that if a point is outside the amoeba of $V$, there is an 
element of the defining ideal which witnesses this 
fact by being lopsided.  This condition is necessary and sufficient
for amoebas of arbitrary codimension, as well as for compactifications
of amoebas inside any toric variety.  Our approach naturally leads to
methods
for approximating hypersurface amoebas and their spines by systems of
linear inequalities.  Finally, we remark
that our main result 
can be seen a precise analogue of a 
Nullstellensatz statement for tropical varieties.
\end{abstract}


\section{Introduction}


\subsection{Statement of results}

Let $V \subset (\CC^*)^r$ be an algebraic variety, defined by an ideal
$I \subset \basering$.  
\begin{definition}[Gel'fand-Kapranov-Zelevinsky \cite{GKZ}]
\rm
The {\bf amoeba} of $V$, is defined to be the image
of $V$ under the map $\Log:(\CC^*)^r \to \RR^r$ defined at the
point $\boldz = (z_1, \ldots, z_r)$ by
$$\Log(\boldz) = (\log|z_1|, \ldots, \log|z_r|).$$
\end{definition}
We denote the amoeba of $V$ by either of $\Am_V$ or $\Am_I$.  If 
$V=Z_f$ is a hypersurface, the zero locus of a single function $f$, we 
will also use the notation $\Am_f$.  We refer the reader to Mikhalkin's
survey article \cite{M} for a broad discussion of amoebas and their 
applications.

Consider $f \in \basering$ and a point $\bolda \in \RR^r$.  Write $f$ as a sum 
of monomials 
$f(\boldz) = m_1(\boldz) + \cdots + m_d(\boldz)$.
Define $f\{\bolda\}$ to be the list of positive real numbers
$$f\{\bolda\} := \big\{|m_1(\Log^{-1}(\bolda))|, 
\ldots, |m_d(\Log^{-1}(\bolda))|\big\}.$$
Note that since the $m_i$ are monomials, this is well defined even 
though $\Log$ is not injective.

\begin{definition}
\rm
\label{def:lopsided}
We say that a list of positive numbers is {\bf lopsided},
if one of the numbers is greater than the sum of all the others. 
\end{definition}

Equivalently, a list of numbers $\{b_1, \ldots, b_d\}$ is
 not lopsided, 
if it is possible to choose complex phases $\phi_i$ ($|\phi_i| = 1$), so
that $\sum \phi_i b_i = 0$.  This follows from the triangle inequality.
We also define 
$$\LopAm_f := \big \{ \bolda \in \RR^r
\ |\ \text{$f\{\bolda\}$ is not lopsided} \big \}.$$

One can easily see that if $\bolda \in \Am_f$ then $f\{\bolda\}$ cannot
be lopsided; in other words $\LopAm_f \supset \Am_f$.
Indeed if $f(z) = 0$ then $m_1(z) + \cdots + m_d(z) = 0$
so it is giving a way to assign complex phases to the list
$\big\{|m_1(z)|, \ldots, |m_d(z)|\big\} = f\{\Log(z)\}$ such that the sum is
$0$.  Thus one should think
of $\LopAm_f$ as a crude approximation to the amoeba $\Am_f$.

\begin{example}
\rm
Suppose 
$f(z_1,z_2) = 1+z_1z_2+z_2^2$, and let $\bolda \in \RR^2$. 
For any complex phases $\phi_1, \phi_2$,
there exist 
$(z_1,z_2) \in \Log^{-1}(\bolda)$ such that
$\phi_1|z_1z_2| = z_1z_2$ and $\phi_2|z_2^2| = z_2^2$.  
Thus $\bolda \in \Am_f$ if and only if $\big\{1, |z_1z_2|, |z_2^2|\big\}$
is non-lopsided; i.e. $\Am_f = \LopAm_f$.
\end{example}

\begin{figure}[htbp]
  \begin{center}
    \epsfig{file=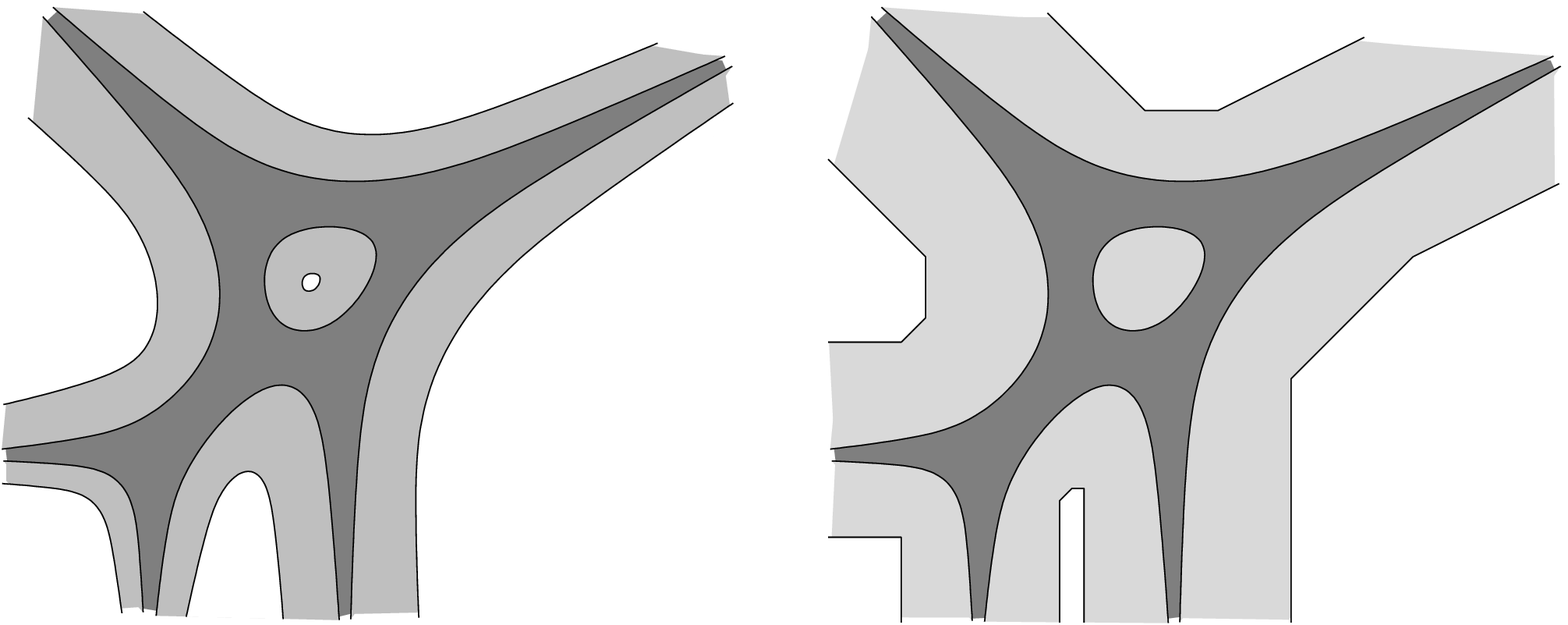,width=6in}
    \caption{The image on the left depicts $\LopAm_f \supset \Am_f$,
     while the right depicts $\SlopAm_f \supset \Am_f$.  Here $\SlopAm_f$
     is not homotopic to $\Am_f$, and in general $\Am_f$ need not be 
     either.}
    \label{fig:la}
  \end{center}
\end{figure}

In the above example
we have enough freedom to choose the phases of the monomials 
$m_i(\boldz)$
for $\boldz \in \Log^{-1}(\bolda)$ that $\LopAm_f = \Am_f$.
However, this only works because $f$ has very few non-zero terms.
In general $\LopAm_f$ can be quite different
from $\Am_f$ (see Figure \ref{fig:la}).  
Nevertheless, we will show that for a suitable multiple of
$f$, we can use this lopsidedness test to get very good approximations
for $\Am_f$.


Let $n$ be a positive integer.  We
consider the polynomials
$$\cyclres{f}(\boldz) =
\prod_{k_1=0}^{n-1} \cdots \prod_{k_r=0}^{n-1}
f(e^{2 \pi i\,k_1/n}\,z_1, \ldots , e^{2 \pi i\,k_r/n}\ z_r).$$
These $\cyclres{f}$ are cyclic resultants 
$$\cyclres{f}(\boldz) =
\Res\Big(\Res\big( \ldots \Res\big(f(u_1z_1, \ldots, u_rz_r), u_1^n-1\big) 
\ldots, u_{r-1}^n-1\big), u_r^n-1\Big)$$
and as such can be practically computed.
Our main result for amoebas
of hypersurfaces is roughly the following.  The precise version is
stated and proved in Section~\ref{subsec:hypersurface}.

\begin{othertheorem}{Theorem~\ref{thm:hypersurface} (rough version)}
As $n \rightarrow \infty$,
the family $\LopAm_{\cyclres{f}}$ converges uniformly to $\Am_f$.
There exists an integer $N$ such that to compute $\Am_f$ to 
within $\varepsilon$, it suffices to compute
$\LopAm_{\cyclres{f}}$ for any $n \geq N$.  
Moreover, $N$ depends only on $\varepsilon$ and the Newton polytope
(or degree) of $f$, and can be computed explicitly from these data.
\end{othertheorem}

This will lead us to the following characterisation of the amoeba 
of a general subvariety of $(\CC^*)^r$.

\begin{othertheorem}{Theorem~\ref{thm:generalamoeba}}
Let $I \subset \basering$ be an ideal.  A point $\bolda \in \RR^r$ is
in the amoeba $\Am_I$ if and only if $g\{\bolda\}$ is not lopsided
for every $g \in I$.
\end{othertheorem}

Phrased another way, if a point $\bolda$ is outside the amoeba $\Am_I$,
a polynomial $f \in I$ may witness this fact by being lopsided
at $\bolda$.  Theorem~\ref{thm:generalamoeba}
then states that there is always a witness.  
We will actually show something slightly stronger, in both 
Theorems~\ref{thm:hypersurface} and~\ref{thm:generalamoeba}.  We will show that 
there is a witness $f$ such that $f\{\bolda\}$ is ``superlopsided'', 
according to the following definition:

\begin{definition}
\rm
Let $d' \geq d \geq 2$.
We say that a list of positive numbers $\{b_1, \ldots, b_{d+1}\}$ is 
\hbox{\bf $d'$-superlopsided} if there exists some $i$ such that 
$b_i > d'b_j$ for all $j \neq i$.  If $d'=d$, we will simply
say the list is superlopsided.
\end{definition}

As before, we also define
$$\SlopAm_f := \big \{ \bolda \in \RR^r
\ |\ \text{$f\{\bolda\}$ is not superlopsided} \big \}.$$

If a list of positive numbers is superlopsided, it is certainly lopsided;  
hence $\SlopAm_f \supset \LopAm_f \supset \Am_f$ (see Figure \ref{fig:la}).
David Speyer observed that each component of the complement of
$\SlopAm_f$ is given by a system of \emph{linear} inequalities, 
making it easier than $\LopAm_f$ to compute explicitly.  Hence 
Theorem~\ref{thm:hypersurface} actually prescribes a method for approximating
$\Am_f$ to within $\varepsilon$ by systems of linear inequalities.  
Similar ideas lead to a method for approximating the {\em spine} of
a hypersurface amoeba.  We discuss these constructions
in Section~\ref{sec:approx}.

The motivation for these results comes from tropical algebraic 
geometry, and from this viewpoint, lopsidedness (rather than
superlopsidedness) is the more natural condition to consider.  
In tropical algebraic geometry, we work with the
semiring $\Rtrop(\odot, \oplus)$.  This a semiring whose underlying 
set is $\RR$, but whose operations are given by 
\begin{itemize}
\item $a \odot b := a + b$
\item $a \oplus b := \max(a,b)$.
\end{itemize}
The operations $\odot$ and $\oplus$ are known as tropical addition
and tropical multiplication.  One can easily check
that they satisfy the usual
commutative, associative, and distributive laws; however there are no
additive inverses. 

A polynomial $g \in \Rtrop[x_1, \ldots, x_r]$ is therefore a piecewise
linear function on $\RR^r$: if 
$$g(\boldx) = 
\bigoplus_{k_1,\ldots, k_r} c_{k_1, \ldots, k_r} 
\odot x_1^{k_1} \odot \cdots \odot x_r^{k_r}$$
then, translated into the usual operations on $\RR$, 
$$g(\boldx) = 
\max \{ c_{k_1, \ldots, k_r} + k_1 x_1 + \cdots + k_r x_r \}.$$
The tropical variety associated to $g$ is then defined to be the singular 
locus of this piecewise linear function.  A tropical variety associated
to a single polynomial $g$ in this way is called a 
{\bf tropical hypersurface}. 

Thus there is a simple Nullstellensatz\footnote{We use the term 
in the literal sense of being a statement about zeros; 
these results are not an analogue of Hilbert's Nullstellensatz.}
for tropical hypersurfaces.
A point $\boldx$ is outside the 
tropical variety of $g$, if there is a single monomial term of $g$ 
which is 
strictly larger than each of the others when evaluated at the point 
$\boldx$.  In terms of the tropical operations this term is strictly
greater 
than the tropical sum of the other terms 
(c.f.  Definition~\ref{def:lopsided}).

More generally the principal results in this paper can be seen of as an
analytic analogue of a theorem of Speyer and Sturmfels \cite{SS} for 
tropical varieties of arbitrary codimension, also known as non-Archimedian 
amoebas.  We discuss this connection in Section~\ref{sec:tropical}.

\subsection{Acknowledgments}
I am deeply grateful to David Speyer, for providing a number of useful
insights, and in particular for drawing my attention to the existence
of Lemma~\ref{lem:averaging}.  I would also like to thank Allen
Knutson providing feedback on early drafts, Matthias Beck for 
answering my questions on Ehrhart polynomials, and Bernd Sturmfels 
for being particularly supportive of this project and for suggesting
the title.  This research was partially supported by NSERC.


\section{The case $r=1$}


\subsection{A heuristic argument}

The idea of the one variable case is simple enough.  Suppose that
$f(z) = \prod_{i=1}^d (z-\alpha_i)$, and for sake of argument, 
assume that the absolute values of the $\alpha_i$ are all distinct,
say $|\alpha_1| > \dots > |\alpha_d| > 0$.  
Then

\begin{alignat*}{2}
\cyclres{f}(z) &= \prod_{i=1}^d && (z^n - \alpha_i^n)  \\
&= z^{nd} &&- (\alpha_1^n + \bdotsb) z^{n(d-1)} \\
&    &&+(\alpha_1^n \alpha_2^n + \bdotsb) z^{n(d-2)} \\
&    &&+\ \ \cdots \ \ +\\
&    &&\pm (\alpha_1^n \cdots \alpha_{d-1}^n + \bdotsb) z^n \\
&    &&\mp \alpha_1^n \cdots \alpha_d^n. \\
\end{alignat*}

For $n$ large, the terms $\bdotsb$ are small in comparison with
the other terms and so this is approximately
$$ g_n(z) = (z^d)^n - (\alpha_1 z^{d-1})^n + (\alpha_1 \alpha_2 z^{d-2})^n
+ \cdots \pm (\alpha_1 \cdots \alpha_{d-1} z)^n
\mp (\alpha_1 \cdots \alpha_d)^n.$$
Suppose $|\alpha_{k+1}| < |z| < |\alpha_k|$. 
Consider $g_n(z) / (\alpha_1 \cdots \alpha_{d-k} z^{d-k})^n$:  as
$n \rightarrow \infty$, every term tends to $0$ except for the 
constant term which is $1$.  Thus for $n$ large, there is a single 
term in $g_n(z)$ and likewise in $\cyclres{f}$ which is much bigger 
in absolute value than all the others.


\subsection{The one variable lemmas}
We now formalise this heuristic argument in a way that will be
useful in proving Theorem~\ref{thm:hypersurface}.
At the crux of the heuristic argument are the following three key facts 
about $\cyclres{f}$:
\begin{enumerate}
\item It has no roots inside a certain annulus (in the heuristic
argument the annulus is
$\{z \in \CC\ \big|\ |\alpha_{k+1}| < |z| < |\alpha_k|$\}).
\item The only non-zero terms which appear are of the form $c z^{nk}$.
\item The number of terms is not too large. (This approach
fails if instead of $\cyclres{f}(z)$, we try to use $\cyclres{f}(z)^D$ 
with $D >> n$.)  
\end{enumerate}

To get a result which we can apply to the multivariable case, we
need to be able to make a uniform statement about polynomials with
these properties.  This is precisely captured by the next two lemmas.  
By applying Lemma~\ref{lem:onevariable} directly to the family of 
functions $\cyclres{f}(z)$
one immediately obtains a complete proof of Theorem~\ref{thm:hypersurface} 
in the $r=1$ case.

\begin{lemma}
\label{lem:goodestimate}
Let $A = \{z\ |\ \beta_0 < |z| < \beta_1)$ be an open annulus in $\CC$. 
We allow the possibilities of $\beta_0 = 0$, and $\beta_1 = \infty$.
Let $\gamma$ be a real number with $\sqrt{\beta_0/\beta_1} \leq \gamma < 1$,
and let $K$ denote the closed interval 
$[\log(\gamma^{-1} \beta_0), \log(\gamma \beta_1)]$.

Suppose $f(z)$ is a polynomial such that 
\begin{enumerate}
\item $f(z)$ has no roots in $A$, and 
\item The only terms which appear in $f(z)$ are of the form $c z^{nl}$,
i.e.  $f(z) = \sum_{l=0}^d m_l(z)$, where $m_l(z) = c_l z^{n(d-l)}$, for some
positive integer $n$.
\end{enumerate}
Then there is some $k$ such that 
\begin{equation}
\label{eqn:goodestimate}
\frac{ |m_l(z_0)|}{|m_k(z_0)|} 
< \frac{\sum_{w \geq |k-l|}  \frac{(d \gamma^n)^w}{w!}}
{2-e^{d \gamma^n}}
\end{equation}
for all $z_0 \in \Log^{-1}(K)$.
\end{lemma}

A similar statement can be made for Laurent 
polynomials: if $f(z)$ is a
Laurent polynomial satisfying properties 1 and 2, and 
$d = \frac{1}{n}(\mindeg(f) - \maxdeg(f))$, the conclusion
\eqref{eqn:goodestimate} remains valid.
Here $\maxdeg(f)$ and $\mindeg(f)$ refer 
respectively to the largest and 
smallest exponents which appear in $f(z)$.  The notation in our proof 
assumes that $f(z)$ is actually a polynomial.

\begin{proof}
We can write
$$f_n(z) = \prod_{i=1}^d (z^n+\alpha_i^n),$$
where $|\alpha_1| \geq \cdots \geq |\alpha_d|$.  We adopt the convention
that $\alpha_0 = 0$, and $\alpha_{d+1} = \infty$.

Since $f_n(z)$ has no roots in $A$ we have 
\begin{equation}
\label{eqn:definingk}
\alpha_{k+1} \leq \beta_0 < \beta_1 \leq \alpha_{k}
\end{equation}
for some $k$ ($0 \leq k \leq d$).  We show that this $k$ has the
desired property.

Consider $g(z) = f_n(z)/ (\alpha_1^n \cdots \alpha_k^n z^{(d-k)n})$.
Write 
$$g(z) = \sum_{l=0}^{d} m'_l(z),$$
where $m'_l(z)$ is a monomial of degree $n(k-l)$.
We now estimate the size of $|m'_l(z_0)|$ at some point 
$z_0 \in \Log^{-1} K$.  If $l=k$, \eqref{eqn:goodestimate} clearly
holds, so we only need to consider $l \neq k$.
First, we will show that if $l \neq k$, we have
\begin{equation}
\label{eqn:lestimate2}
|m'_l(z_0)| < \sum_{w \in \ZZ} 
{k \choose l-w}{d-k \choose w} \gamma^{n(2w + k-l)}.
\end{equation}

Certainly for all $l \neq k$ we have
\begin{align}
|m'_l(z_0)| &=
\frac{
\left | \left (
\sum_{s_1 < \cdots <s_l} \alpha_{s_1}^n \cdots \alpha_{s_l}^n
\right ) z^{n(d-l)} \right |
}
{ |(\alpha_1^n \cdots \alpha_k^n z^{(d-k)n} | } \notag\\
& \leq \sum_{s_1 < \cdots < s_l} 
\left |
\frac{\alpha_{s_1}^n \cdots \alpha_{s_l}^n}
{\alpha_1^n \cdots \alpha_k^n}  z^{(k-l)n}
\right |. \label{eqn:lestimate1}
\end{align}
We use the following facts:
\begin{itemize}
\item If $s_i \leq k$, we have $|\alpha_{s_i}/\alpha_i| \leq 1$
(since $s_i \geq i$).
\item If $s_i > k$ and $i\leq k$, we have 
$|\alpha_{s_i}/\alpha_i| \leq |\beta_0/\beta_1| \leq \gamma^2$.
\item If $i \leq k$, then $|z_0/\alpha_i| \leq |z_0/\beta_1| \leq \gamma$.
\item If $i > k$, then $|\alpha_i/z_0| \leq |\beta_0/z_0| \leq \gamma$.
\end{itemize}
If $l<k$ we write \eqref{eqn:lestimate1} as
\begin{align*}
|m'_l(z_0)| &\leq \sum_{s_1 < \cdots < s_l}
\left | \frac{\alpha_{s_1}}{\alpha_1} \right |^n
\cdots
\left | \frac{\alpha_{s_l}}{\alpha_l} \right |^n
\left | \frac{z_0}{\alpha_{l+1}} \right |^n
\cdots
\left | \frac{z_0}{\alpha_k} \right |^n \\
&\leq \sum_{s_1 < \cdots < s_l}
(1^n)^{\# \{i\ |\ s_i \leq k\}} (\gamma^{2n})^{\# \{i\ |\ s_i > k\}} (\gamma^n)^{k-l}
\\
&= \sum_{w \in \ZZ} 
{k \choose l-w}{d-k \choose w} \gamma^{n(2w + k-l)}.
\end{align*}
This last equation comes from the fact that 
${k \choose l-w}{d-k \choose w}$ is
the number of
$s_1 < \cdots < s_l$ with $\# \{i\ |\ s_i \leq k\} = w$.
Similarly, if $l>k$ we write
\begin{align*}
|m'_l(z_0)| &\leq \sum_{s_1 < \cdots < s_l}
\left | \frac{\alpha_{s_1}}{\alpha_1} \right |^n
\cdots
\left | \frac{\alpha_{s_k}}{\alpha_k} \right |^n
\left | \frac{\alpha_{k+1}}{z_0} \right |^n
\cdots
\left | \frac{\alpha_l}{z_0} \right |^n \\
&\leq \sum_{s_1 < \cdots < s_l}
(1^n)^{\# \{i\ |\ s_i \leq k\}} 
(\gamma^{2n})^{\# \{i\ |\ s_i > k \hbox{ and } i \leq k\}} (\gamma^n)^{l-k}
\\
&= \sum_{w \in \ZZ} 
{k \choose l-w}{d-k \choose w} \gamma^{n(2(w+k-l)+l-k}
\\
&= \sum_{w \in \ZZ} 
{k \choose l-w}{d-k \choose w} \gamma^{n(2w + k-l)},
\end{align*}
which establishes \eqref{eqn:lestimate2}.

Thus we obtain 
\begin{align*}
\sum_{l \neq k} |m'_l(z_0)| 
& \leq \sum_{l \neq k} \sum_{w \in \ZZ} 
{k \choose l-w}{d-k \choose w} \gamma^{n(2w + k-l)}  \\
& = \sum_{l \neq k} \sum_{w \in \ZZ} 
{k \choose w+k-l}{d-k \choose w} \gamma^{n(2w + k-l)} \\
& = \sum_{w \in \ZZ} \sum_{w' \neq w}
{k \choose w'}{d-k \choose w} \gamma^{n(w+w')} \\
& \leq \sum_{(w,w') \in \ZZ^2 \setminus (0,0)}
{k \choose w'}{d-k \choose w} \gamma^{n(w+w')} \\
& = -1 + \sum_{w'=0}^k {k \choose w'} \gamma^{nw'}
\sum_{w=0}^{d-k} {d-k \choose w} \gamma^{nw} \\
& = -1 + (1+\gamma^n)^k (1+\gamma^n)^{d-k} \\
& = (1+\gamma^n)^d -1.
\end{align*}
In particular for each $l \neq k$ we have
\begin{equation}
\label{eqn:badestimate}
|m'_l(z_0)| \leq (1+\gamma^n)^d -1.
\end{equation}

However, we can do 
slightly better than this, by noting that
the smallest power of $\gamma$ that appears on the right side of the
inequality \eqref{eqn:lestimate2}
is $\gamma^{|k-l|}$.  Thus, whereas \eqref{eqn:badestimate} tells us
that $|m'_l(z_0)| < \sum_{w \geq 1} {d \choose w} \gamma^{nw}$, in fact
we have:
\begin{align}
|m'_l(z_0)| &< \sum_{w \geq |k-l|} {d \choose w} \gamma^{nw} \notag \\
&< \sum_{w \geq |k-l|}  \frac{(d \gamma^n)^w}{w!}. 
\label{eqn:lestimate3}
\end{align}
(Although \eqref{eqn:lestimate2} is a better estimate than 
\eqref{eqn:lestimate3}, the latter will prove to be more useful to us.)

For $m_k(z_0)$ we have the following estimate:
\begin{align}
|m'_k(z_0)| &=
\frac{
\left | \left (
\sum_{s_1 < \cdots <s_k} \alpha_{s_1}^n \cdots \alpha_{s_k}^n
\right ) z^{n(d-k)} \right |
}
{ |(\alpha_1^n \cdots \alpha_k^n z^{(d-k)n} | } \notag \\
& \geq 1 - \sum_{\substack{s_1 < \cdots < s_k\\s_k>k}}
\left | \frac{\alpha_{s_1}}{\alpha_1} \right |^n
\cdots
\left | \frac{\alpha_{s_k}}{\alpha_k} \right |^n \notag \\
& \geq 1 - \sum_{\substack{s_1 < \cdots < s_k\\s_k>k}}
(1^n)^{\# \{i\ |\ s_i \leq k\}} (\gamma^{2n})^{\# \{i\ |\ s_i > k\}}
\notag \\
&= 1 - \sum_{w \geq 1} 
{k \choose k-w}{d-k \choose w} \gamma^{2nw} \notag \\
&= 1 - \sum_{w \geq 1} 
{k \choose w}{d-k \choose w} \gamma^{2nw} \notag \\
& \geq 1 - \sum_{(w,w') \in \ZZ^2 \setminus (0,0)}
{k \choose w'}{d-k \choose w} \gamma^{n(w+w')} \notag \\
& \geq 1 - \big ( (1+\gamma^n)^d -1 \big ) \label{eqn:kestimate1} \\
& \geq 2-e^{d \gamma^n}. \label{eqn:kestimate2}
\end{align}

Combining the estimates \eqref{eqn:lestimate3} and 
\eqref{eqn:kestimate2} we obtain
\begin{align*}
\frac{|m_l(z_0)|}{|m_k(z_0)|} &= 
\frac{|m'_l(z_0)|}{|m'_k(z_0)|}  \\
&< \frac{\sum_{w \geq |k-l|}  \frac{(d \gamma^n)^w}{w!}}
{2-e^{d \gamma^n}}
\end{align*}
for $l \neq k$, as required.  

\end{proof}


\begin{remark}
\label{rmk:domterm}
\rm
It should be noted that we have actually determined which term is the
special term $m_k$.  This is done in \eqref{eqn:definingk}.  
If $f_n$ is a polynomial, then $n(d-k)$ is the number of roots (counted with 
multiplicity) of $f_n$ inside the disc $\{ |z| \leq \beta_0 \}$.
If $f_n$ is a Laurent polynomial, $n(d-k) - \mindeg(f_n)$ will be the number 
of roots inside $\{0 < |z| \leq \beta_0\}$.
\end{remark}

\begin{lemma}
\label{lem:onevariable}
As before, let $A = \{z\ |\ \beta_0 < |z| < \beta_1)$.  Let
$\gamma$ be a real number with $\sqrt{\beta_0/\beta_1} \leq \gamma < 1$,
and $K = [\log(\gamma^{-1} \beta_0), \log(\gamma \beta_1)]$.
Fix positive integers $c_0, D_0, c_1, D_1$.

Suppose $f_n(z)$ is a family of polynomials such that
\begin{enumerate}
\item $f_n(z)$ has no roots in $A$,
\item The only terms which appear in $f_n(z)$ are of the form $c z^{nk}$.
\item $\deg(f_n) \leq c_0 n^{D_0+1}$.  More generally, if $f_n$ are
Laurent polynomials, we can ask that 
$\maxdeg(f_n) - \mindeg(f_n) \leq c_o n^{D_0+1}$.
\end{enumerate}
Assume $n$ is large enough that
$n \log \gamma^{-1} \geq (D_0+D_1) \log n + \log(8/3 c_0c_1)$.  
Then $f_n\{a\}$ is $(\mdeg{1})$-superlopsided for all $a \in K$.
\end{lemma}

%
%

\begin{proof}
From Lemma~\ref{lem:goodestimate}, if we write 
$f_n(z) = \sum_{l=0}^d m_l(z)$, there is some $k$ 
such that for $l \neq k$,
$$ \frac{|m_l(z_0)|}{|m_k(z_0)|} 
\leq \frac{(1+\gamma^n)^d-1}{2 - (1+\gamma^n)^d}.$$
Here we are using the weaker estimate in \eqref{eqn:badestimate}
to bound the numerator, and \eqref{eqn:kestimate1} for the denominator.  
Thus,
\begin{align*}
\mdeg{1}\frac{|m_l(z_0)|}{|m_k(z_0)|} & \leq 
\frac{\mdeg{1} \big( (1+\gamma^n)^d-1 \big)}{2 - (1+\gamma^n)^d}  \\
& \leq \frac{\mdeg{1} \big( (1+\gamma^n)^{\mdeg{0}} - 1 \big)}
{2 - (1+\gamma^n)^{\mdeg{0}}} \\
& \leq \frac{1/2}{2-(3/2)} = 1.
\end{align*}
The last step uses Calculation~\ref{cal:1}, which can be found
in Appendix~\ref{app:calcdetails}.
This establishes that $f_n\{\Log(z_0)\}$ is $(\mdeg{1})$-superlopsided.
\end{proof}


\subsection{Accuracy of bounds}

This proof of Lemma~\ref{lem:onevariable} 
uses a number of inequalities to give an answer to the 
question `how large does $n$ have to be?'  Many of these inequalities, at 
first glance, appear not to be very tight.  The answer we obtain is 
certainly
not the best possible, and in most situations we should not expect to
need quite such a large $n$.  However, without assuming more about $f_n$ it 
is surprisingly close to the best answer, particularly as 
$\gamma \rightarrow 1$.
The best general answer for the question will be of the same 
form, i.e. 
$$n \log \gamma^{-1} \leq (D_0 + D_1) \log n + \log \Theta(c_0c_1).$$
One can see this by performing the requisite analysis on
the polynomials $h_n(z) = (z^n + 1)^{\mdeg{0}}$.  
As a general heuristic, the more
closely the roots of $f_n$ are packed, the larger $n$ will have to be; 
thus the family of polynomials $h_n(z)$, where
every root has as high a multiplicity as possible, is where we should
expect our worst case behaviour to occur.  

Suppose we want $n$ large enough to guarantee that 
$h_n\{\log |z_0|\}$ is 
$(\mdeg{1})$-super\-lopsided for $|z_0|<\gamma < 1$.  Write
$h_n(z) = 1 + \mdeg{0} z^n + \cdots$.  We know $1$ will be the dominant
term as $n$ gets large (since $1 = \lim_{n \to \infty} h_n(z_0)$), thus 
we need $(\mdeg{1})(\mdeg{0} z_0^n) < 1$, or equivalently
$$n \log \gamma^{-1} \geq (D_0+D_1) \log n + \log (c_0 c_1).$$

If we only want to guarantee that $f_n\{a\}$ is lopsided for $a \in K$, we 
need $n$ large enough so that
$$\sum_{l \neq k} |m_l(z_0)| \leq (1+ \gamma^n)^d - 1 <
2 - (1 + \gamma^n)^d \leq |m_k(z_0)|,$$  
or equivalently $(1+ \gamma^n)^{\mdeg{0}} \leq 3/2.$
This will hold if we have
$$n \log \gamma^{-1} \geq D_0 \log n + \log\left(\frac{c_0}{\log 3/2}\right).$$
So $n$ only
needs to be about half as big to guarantee that $f_n\{a\}$ is lopsided
as it does to guarantee that $f_n\{a\}$ is superlopsided.  
Again, we can show see that this is
fairly close to the best answer by considering $(z^n + 1)^{\mdeg{0}}$.


\section{The hypersurface case}


\subsection{Preliminaries}
\label{subsec:prehypersurface}

In this section we prove our main theorem characterising the amoeba
of a hypersurface.  If $f(\boldz) \in \basering$, we consider the 
Laurent polynomials 
$$\cyclres{f}(\boldz) =
\prod_{k_1=0}^{n-1} \cdots \prod_{k_r=0}^{n-1}
f(e^{2 \pi i\ k_1/n}\ z_1, \ldots , e^{2 \pi i\ k_r/n}\ z_r).$$
Theorem~\ref{thm:hypersurface} states that for $\varepsilon>0$, and
a point $\bolda \in \RR^r$ in the complement of the amoeba $\Am_f$ whose 
distance from $\Am_f$ is at least $\varepsilon$, we can choose
$n$ large enough so that $\cyclres{f}\{\bolda\}$ is superlopsided.
Moreover, the theorem gives a upper bound on how large
$n$ needs to be,
based only on $\varepsilon$ and the {\em Newton polytope} of $f$.

The idea behind the proof of Theorem~\ref{thm:hypersurface}
is to look at the family of $\cyclres{f}(\boldz)$
and interpret this as a function of
a single variable $z_i$.
At the point $\boldzeta = (\zeta_1, \ldots, \zeta_r) \in \CC^r$, we define
$$\projcyclres{f}{i}(z) := \cyclres{f}(\zeta_1, \ldots, \zeta_{i-1}, z,
\zeta_{i+1}, \ldots, \zeta_n).$$
We will apply Lemma~\ref{lem:onevariable} to 
these, and a find a single dominant term in this polynomial of
one variable.  Then by an averaging argument, we show that this implies 
that $\cyclres{f}$ has a single dominant term.  

First, however, we need a few simple observations.

\begin{proposition}
\label{prop:sameamoeba}
$\Am_f = \Am_{\cyclres{f}}.$
\end{proposition}

\begin{proof}
$\cyclres{f}(\boldz)$ is a product of terms
$g_{u_1, \ldots, u_r}(\boldz) = f(u_1 z_1, \ldots, u_r z_r)$, 
where $u_i^n = 1$. Since $|u_i|=1$, $\Am_{g_{u_1, \ldots , u_r}} = \Am_f$, 
and so 
$\Am_{\cyclres{f}} = \bigcup \Am_{g_{u_1, \ldots , u_r}} = \Am_f$.
\end{proof}

We will also need to know some information about the number and degree
of the terms which appear in $\cyclres{f}$.  First note the following
important fact.

\begin{proposition}
\label{prop:powersofn}
The only monomials which appear in $\cyclres{f}$ are of the form
$c\,z_1^{nk_1}\cdots z_r^{nk_r}$.
In particular the only terms appearing in 
$\projcyclres{f}{i}(z)$ are of the form $c z^{nk}$.
\end{proposition}

\begin{proof}
Let $C_n$ denote the cyclic group of roots of $z^n-1$.
$\cyclres{f}$ is manifestly invariant under the group action of
$(C_n)^r$ acting on $\basering$ by 
$(u_1, \ldots, u_n) \cdot g(\boldz) = g(u_1 z_1, \dots, u_n z_n)$.
Thus each monomial of $\cyclres{f}$ must be invariant under this action.
The only monomials with this property are of the form
$c\,z_1^{nk_1}\cdots z_r^{nk_r}$.

The statement about $\projcyclres{f}{i}(z)$ follows immediately.
\end{proof}

Recall that if $g \in 
\basering$ its {\bf Newton polytope}, denoted $\Delta(g)$,
is the subset of $\RR^r$ defined as the convex hull of the exponent 
vectors of the monomials which appear in $g$.

For any polytope $\Delta$, let $d(\Delta)$ be any upper bound on 
$(\# \{\ZZ^r \cap m \Delta\})/m^r$.  In general, it is not easy
to find a tight upper bound for this number.  If one can compute
the Ehrhart polynomial of $\Delta$ explicitly, than an easy upper
bound is the sum of the positive coefficients.  Otherwise, it is
possible to bound the coefficients of the Ehrhart polynomial in
terms of the volume of $\Delta$ (see \cite{BM}).  Using these
estimates, for each $r$ one can compute constants $A$ and $B$ such that
$(\# \{\ZZ^r \cap m \Delta\})/m^r < A \cdot \Vol(\Delta) + B$.

Clearly we have $\Delta(\cyclres{f}) = n^r \Delta(f)$.  
This gives us an upper bound on the number of terms that $\cyclres{f}$
can have:

\begin{proposition}
\label{prop:newtonpolytope}
Let $d = d(\Delta(f))$. Then $\cyclres{f}$ has at most $d\,n^{r^2-r}$ terms.
\end{proposition}

\begin{proof}
By Proposition~\ref{prop:powersofn}, the number of terms in $\cyclres{f}$ 
is at most the number of integral points in 
$\frac{1}{n}\Delta(\cyclres{f}) = n^{r-1}\Delta(f)$. This
is less than or equal to $d\,n^{r^2-r}$.
\end{proof}

Finally we need to know something about 
$\maxdeg(\projcyclres{f}{i}) - \mindeg(\projcyclres{f}{i})$.
Let $$c_i(f) := \max x_i(\Delta(f)) - \min x_i(\Delta(f)),$$ 
where $x_i$ denotes the $i^{\rm{th}}$ coordinate function on $\RR^r$.

\begin{proposition}
\label{prop:projnewtonpolytope}
$\maxdeg(\projcyclres{f}{i}) - \mindeg(\projcyclres{f}{i}) = c_i(f) n^r.$
\end{proposition}

\begin{proof}
\begin{align*}
\maxdeg(\projcyclres{f}{i}) - \mindeg(\projcyclres{f}{i}) 
&= \max x_i(\Delta(\cyclres{f})) - \min x_i(\Delta(\cyclres{f}))  \\
&= n^r \max x_i(\Delta(\cyclres{f}))  - n^r \min x_i(\Delta(\cyclres{f})) \\
&= d_i n^r.
\end{align*}
\end{proof}


\subsection{Proof of Theorem~\ref{thm:hypersurface}}
\label{subsec:hypersurface}

Armed with these facts and Lemmas~\ref{lem:goodestimate} 
and~\ref{lem:onevariable} we are
now in a position to precisely state and prove our main result for 
amoebas of hypersurfaces.

\begin{theorem}
\label{thm:hypersurface}
Let $\varepsilon>0$.
Suppose $\bolda = (a_1, \ldots, a_r) \in \RR^r \setminus \Am_f$ is a point 
in the amoeba complement whose distance from $\Am_f$ is at least 
$\varepsilon$.  
Let 
$d = d(\Delta(f))$, and
$c = \max \{c_i(f)\ |\ 1 \leq i \leq r\}$.
\begin{enumerate}
\item
If $n$ is large enough so that
\begin{equation}
\label{eqn:nlarge1}
n \varepsilon \geq (r-1) \log n + \log \big((r+3)2^{r+1}c\big)
\end{equation}
then $\cyclres{f}\{\bolda\}$ is lopsided.
\item
If $n$ is large enough so that
\begin{equation}
\label{eqn:nlarge2}
n \varepsilon \geq (r^2-1) \log n + \log(16/3 cd)
\end{equation}
then $\cyclres{f}\{\bolda\}$ is superlopsided.  (In fact it is
$(d n^{r^2-r})$-superlopsided).
\end{enumerate}
\end{theorem}

The key to reducing to the one variable case is the following basic
result from complex analysis.
\begin{lemma}
\label{lem:averaging}
Let $f(\boldz)$ be a Laurent polynomial, and write 
$f(\boldz) = \sum_\vecj m_\vecj(\boldz)$, where 
$m_\vecj(\boldz) = m_{j_1, \ldots, j_r}(\boldz) =
b_{j_1, \ldots, j_r} z_1^{j_1} \cdots z_r^{j_r}$. 
Suppose for all $\boldzeta \in \Log^{-1}(\bolda)$
we have $|f(\boldzeta)| \leq M$.  Then for each $\vecl$,
$ |m_\vecl(\boldzeta)| \leq M$.
\end{lemma}

\begin{proof}
We integrate the equations $M \geq |f(\boldzeta)|$ over the set
$\Log^{-1}(a_1, \ldots, a_r)$.

\begin{align*}
M & \geq 
\frac{1}{(2\pi)^r}
\int_{\theta_1=0}^{2\pi} \cdots \int_{\theta_r=0}^{2\pi} 
\left|
\sum_\vecj m_\vecj (e^{a_1+i\theta_2}, \ldots, e^{a_r+i\theta_r})
\right|\ d\theta_1 \cdots d\theta_r \\
& \geq
\left|
\frac{1}{(2\pi i)^r}
\int_{|z_1|=1} \cdots \int_{|z_r|=1}
\sum_\vecj 
\frac{m_\vecj (e^{a_1}z_1, \ldots, e^{a_r}z_r)}
{z_1^{l_1} \cdots z_r^{l_r}}
\frac{d z_1}{z_1} \cdots \frac {d z_1}{z_1}
\right| \\
& = | m_\vecl (e^{a_1}, \ldots, e^{a_r}) | \\
& = | m_\vecl (\boldzeta) | .
\end{align*}

\end{proof}

%
%

\begin{proof}[Proof of Theorem~\ref{thm:hypersurface}]
Let $\gamma = e^{-\varepsilon}$, and
let $A_i = \{z\ |\ \gamma e^{a_i} < z < \gamma^{-1} e^{a_i} \}$.

Take $\boldzeta \in \Log^{-1}(\bolda)$.  
The function $\projcyclres{f}{i}$ has no roots in the annulus $A_i$.  If it
did, say at $z_0$, $\cyclres{f}$ would have a corresponding root at
$\boldzeta'=(\zeta_1, \ldots, \zeta_{i-1}, z_0, \zeta_{i+1}, \ldots, \zeta_r)$.
So $\boldzeta' \in \Am_{\cyclres{f}} = \Am_f$, hence by assumption,
$\| \Log(\boldzeta') - \bolda\| > \varepsilon$.  But, 
since $z_0 \in A_i$,
$\| \Log(\boldzeta') - \bolda \| = |\log(z_0) - a_i| < \varepsilon$.

We'll first prove Statement 2.  Assume that \eqref{eqn:nlarge2}
holds.

By Propositions~\ref{prop:powersofn} and~\ref{prop:projnewtonpolytope}, 
we are in a position to apply
Lemma~\ref{lem:onevariable} to the Laurent polynomial
$\projcyclres{f}{i}(z)$, taking $A$, and $\gamma$ as above,
and $c_0 =c$, $D_0 = r-1$, $c_1 = 2d$ and $D_1 = r^2 -r$.
With these choices, $n$ is large enough so that 
\begin{equation}
\label{eqn:issuper}
\text{$\projcyclres{f}{i}\{a_i\}$ is $(2dn^{r^2-r})$-superlopsided.}
\end{equation}
This is true for any $\boldzeta \in \Log^{-1}(\bolda)$, and since
$\Log^{-1}(\bolda)$ is connected,
we see that for any fixed $i$, it must be the same monomial term
which dominates in each $\projcyclres{f}{i}$,
independent of the choice of $\boldzeta$.  Let this be the
$z^{nk_i}$-term, and let $\veck = (k_1, \ldots, k_r)$.

Write 
$$\cyclres{f}(\boldz) =
\sum_\vecj m_\vecj (\boldz),$$
where $m_\vecj(\boldz)$ is the monomial
$b_\vecj z_1^{nj_1}\cdots z_r^{nj_r}$.

Let $M = |m_\veck (\boldzeta)|$.  Note this does not depend
on the particular choice of $\boldzeta$.  Let
$$ \mu = 
\max \big \{|m_\vecl(\boldzeta)|\ 
\big |\ \vecl \neq \veck \}.$$
We wish to show that $\mu< \frac{1}{dn^{r^2-r}}M$.  Since the
number of terms in $\cyclres{f}$ is at most $dn^{r^2-r}$, this will
imply
that $\cyclres{f}$ is superlopsided.

We consider the terms of $\projcyclres{f}{i}(z)$, and
for simplicity of notation, we will temporarily assume that $i=1$.
The $z^{nl}$-term in $\projcyclres{f}{1}(z)$ is
$$\sum_{j_2} \cdots \sum_{j_r} m_{l, j_2, \ldots, j_r}
(z, \zeta_2, \ldots, \zeta_r).$$
For $l \neq k_1$, \eqref{eqn:issuper} tells us that
\begin{align*}
\left|
\sum_{j_2} \cdots \sum_{j_r} m_{l, j_2, \ldots, j_r} (\boldzeta)
\right| 
& < \frac{1}{2 dn^{r^2-r}}
\left| 
\sum_{j_2} \cdots \sum_{j_r} m_{k, j_2, \ldots, j_r} (\boldzeta)
\right|  \\
& \leq \frac{1}{2 dn^{r^2-r}}
\sum_{j_2} \cdots \sum_{j_r} |m_{k, j_2, \ldots, j_r} (\boldzeta)| \\
& \leq \frac{M+ \mu dn^{r^2-r}}{2 dn^{r^2-r}}.
\end{align*}
Since this is true for all $\boldzeta$, by Lemma~\ref{lem:averaging},
$| m_{l, j_2, \ldots, j_r} (\boldzeta) |
< (M+ \mu dn^{r^2-r})/(2 dn^{r^2-r})$
provided $l \neq k_1$.  But, this argument is valid for any $i$, so in 
fact we have
$$ |m_\vecl (\boldzeta) | 
< \frac{M+ \mu dn^{r^2-r}}{2 dn^{r^2-r}}$$
for all $\vecl \neq \veck$.  In particular, since 
$\mu = |m_\vecl (\boldzeta)|$ for some $\vecl \neq \veck$, we deduce
\begin{alignat*}{2}
&& \mu &< \frac{M +\mu dn^{r^2-r}}{2 dn^{r^2-r}} \\
&\Rightarrow\qquad & 2 dn^{r^2-r} \mu &< M + dn^{r^2-r} \mu \\
&\Rightarrow\qquad & \frac{\mu}{M} &< \frac{1}{dn^{r^2-r}} \\
\end{alignat*}
as required.

We now prove Statement 1.  Although the approach is essentially the
same, it is slightly more difficult, and hence requires some additional
lemmas (Calculations~\ref{cal:lopbounds} and~\ref{cal:expderiv}, 
which can be found in Appendix~\ref{app:calcdetails}).  
The reason for this is that
we cannot get these bounds by appealing directly to 
Lemma~\ref{lem:onevariable}.  Instead, we use
Lemma~\ref{lem:goodestimate}, which gives better estimates for
the coefficients 
of $\projcyclres{f}{i}$.

Write $\projcyclres{f}{i}(z) = \sum_j m^i_j(z)$, where $m^i_j(z)=b_j z^{nj}$.
Then for each $i$, 
$\maxdeg \projcyclres{f}{i}(z) - \mindeg \projcyclres{f}{i}(z)
\leq c n^r$.  So by Lemma~\ref{lem:goodestimate} there is some 
$k_i$ such that
$$
\frac{ |m^i_l(\boldzeta)|}{|m^i_{k_i}(\boldzeta)|} 
< \frac{\sum_{w \geq |k_i-l|}  \frac{(c n^{r-1} \gamma^n)^w}{w!}}
{2-e^{c n^{r-1} \gamma^n}} ,
$$
and in fact it is the same $k_i$ for all choices of $\boldzeta$.

As before, let $M = |m_\veck (\boldzeta)|$, and let
$\sigma = \sum_{\vecj \neq \veck} |m_\vecj(\boldzeta)|$.  We have
for any $\boldzeta$ that $|m^i_{k_i}(\boldzeta)| < M + \sigma$; thus
$$
\frac{ |m^i_l(\boldzeta)|}{M+\sigma} 
< \frac{\sum_{w \geq |k_i-l|}  \frac{(c n^{r-1} \gamma^n)^w}{w!}}
{2-e^{c n^{r-1} \gamma^n}}.
$$
This is true for all $\boldzeta \in \Log^{-1}(\bolda)$, thus by 
Lemma~\ref{lem:averaging} we actually have
\begin{equation}
\label{eqn:a}
\frac{|m_\vecl(\boldzeta)|}{M+\sigma}
< \frac{\sum_{w \geq |k_i-l_i|}  \frac{(c n^{r-1} \gamma^n)^w}{w!}}
{2-e^{c n^{r-1} \gamma^n}}.
\end{equation}
The inequality \eqref{eqn:a} is valid for all $i$. In particular we should
take $i$ to be the value which maximises $|k_i - l_i|$, as this will
be the one which minimises the right hand side.

As $\sigma$ is a sum of terms $|m_\vecl(\boldzeta)|$,
can use \eqref{eqn:a} to get an upper bound on $\sigma/(M+\sigma)$.
To do this, note that the right hand side of \eqref{eqn:a} depends
on $\veck$ and $\vecl$ only insomuch as it depends on the maximum
value of $|k_i - l_i|$.  Thus we need to count, for each 
$w_0 \in \ZZ_+$, how many
values of $\vecl$ in this sum satisfy 
$$\max \big \{ |k_i - l_i|\ \big |\ i \in \{1, \ldots, r\} \big \} = w_0.$$
An upper bound for this answer, is the number of integral points in
$(j_1, \ldots, j_r) \in \ZZ^r$ with $\max |j_i| = w_0$.  There are
$(2w_0+1)^r - (2w_0-1)^r$ of these, and this number is in turn less
than $2^r{w_0 + r \choose r}$.  

Using $2^r{w_0 + r \choose r}$
as our upper bound for the number of $\vecl$ with 
$\max \{ |k_i-l_i|\} = w_0$, \eqref{eqn:a} gives us
$$
\frac{\sigma}{M+\sigma} <
\sum_{w_0 \geq 1} 2^r{w_0 + r \choose r}
\frac{ \sum_{w \geq w_0 }
\frac{(c n^{r-1} \gamma^n)^w}{w!}}
{2-e^{c n^{r-1} \gamma^n}},
$$
and by Calculation~\ref{cal:expderiv} (see Appendix~\ref{app:calcdetails})
this becomes
\begin{equation}
\label{eqn:b}
\frac{\sigma}{M+\sigma} <
2^r \left( 
\frac{e^{(r+2) c n^{r-1} \gamma^n} -1}
{2-e^{c n^{r-1} \gamma^n}} \right).
\end{equation}

Assume now that \eqref{eqn:nlarge1} holds.
By Calculation~\ref{cal:lopbounds} (see Appendix~\ref{app:calcdetails}), 
$n$ is large enough so that the right hand
side is less than $1/2$.  Thus we have $\sigma < (M+\sigma)/2$, or 
equivalently, $\sigma < M$, as required.

\end{proof}


\subsection{Corollaries}

As an immediate consequence of Theorem~\ref{thm:hypersurface} we get:

\begin{corollary}
\label{cor:uniform}
Both the families $\{\SlopAm_{\cyclres{f}}\ |\ n \in \ZZ_+\}$
and $\{\LopAm_{\cyclres{f}}\ |\ n \in \ZZ_+\}$
converge uniformly to $\Am_f$.
\end{corollary}

\begin{proof}
We have
$\SlopAm_{\cyclres{f}} \supset 
\LopAm_{\cyclres{f}} \supset \Am_{\cyclres{f}} = \Am_f$, and for
$n$ sufficiently large, the distance from any point
in $\SlopAm_{\cyclres{f}}$ to $\Am_f$ will be less than 
$\frac{1}{n}((r^2-1) \log n + \log(16/3 cd))$, which tends to $0$
as $n \rightarrow \infty$.
\end{proof}

Although Theorem~\ref{thm:hypersurface} only showed that
$\cyclres{f}\{\bolda\}$ is {\em superlopsided} 
for $n$ sufficiently large.
we can easily deduce something stronger.

\begin{corollary}
\label{cor:verylopsided}
Let $\alpha > d$.  If $n$ sufficiently large, then
$\cyclres{f}\{\bolda\}$ is $(\alpha n^{r^2-r})$-superlopsided.
\end{corollary}

Equivalently, if we write $\cyclres{f}\{\bolda\} = 
\{b_{n1}, \ldots, b_{nd_n}\}$, 
then for $n$ sufficiently large, some $b_{nk_n}$ will dominate, and in fact
$$\lim_{n \to \infty} \frac{\sum_{j \neq k_n} b_{nj}}{b_{nk_n}} = 0.$$

\begin{proof}
To see this, one simply has to look at \eqref{eqn:b}.  
As $n \rightarrow \infty$, the right hand side approaches $0$, and thus
so does $\sigma/M$.  

Alternatively, one can note that any $\alpha >d$ can be used in place
of $d$ in Theorem~\ref{thm:hypersurface}.  Thus we see that it 
suffices to take $n$ so that
$$n\varepsilon \geq (r^2-1) \log n + \log(16/3 c\alpha).$$
\end{proof}


\subsection{Accuracy of bounds}

Just as was the case in Lemma~\ref{lem:onevariable}, 
the bounds on $n$ given in Theorem~\ref{thm:hypersurface} are 
not quite optimal: there are a number of places in which the
inequalities can obviously be made tighter.  However, as 
$\varepsilon \rightarrow 0$, the bounds are at least 
asymptotically correct.  

To see this for the superlopsided case, we can consider the example 
$$f(\boldz) = (1-z_1)^{D_1} \cdots (1-z_r)^{D_r}.$$  The amoeba
$\Am_f$ is the union of all coordinate hyperplanes in $\RR^r$.
We can
easily compute
\begin{align*}
\cyclres{f}(\boldz) 
&= (1-z_1^n)^{D_1n^{r-1}} \cdots (1-z_r^n)^{D_rn^{r-1}} \\
&= 1 - D_1n^{r-1}z_1^n - \cdots - D_r n^{r-1}z_r^n + \cdots \ .
\end{align*}
If our point is $\bolda = (a_1, \ldots, a_r)$, with each 
$-\varepsilon < a_i < 0$, 
the dominant term of $\cyclres{f}$, will be $1$.
The number of non-zero monomial terms in $\cyclres{f}$ is greater than 
$D_1 \cdots D_r n^{r^2-r}$; 
Thus to have $\cyclres{f}\{a\}$ superlopsided, we certainly need
\begin{alignat*}{2}
&& D_i n^{r-1} e^{a_in} &< \frac{1}{D_1 \cdots D_r n^{r^2-r}} \\
&\Rightarrow\qquad & e^{-a_in} &> D_i (D_1 \cdots D_r) n^{r^2-r} \\
\end{alignat*}
for all $i$.  This implies that
$$n \varepsilon > (r^2-1) \log n + 
\log \big((D_1 \cdots D_r) \max \{D_i\} \big ).$$
In contrast, \eqref{eqn:nlarge2} says that in this
example we should take $n$ so that
$$n \varepsilon > (r^2-1) \log n +
\log \big(16/3 (D_1+1) \cdots (D_r+1) \max \{D_i\} \big ).$$

Our bound \eqref{eqn:nlarge1} for lopsidedness appears to be
slightly less satisfactory.  In the above example, to guarantee
lopsidedness, one needs $n$ large enough so that 
$$
(1-z_1^n)^{D_1n^{r-1}} \cdots (1-z_r^n)^{D_rn^{r-1}} - 1 < 1.$$
This will hold when
\begin{alignat*}{2}
&& (e{a_1n}+1)^{D_1n^{r-1}} \cdots (e^{a_rn}+1)^{D_rn^{r-1}} &< 2 \\
&\Leftrightarrow\qquad & D_1 n^{r-1} \log(1+e^{a_1n}) +
\cdots D_r n^{r-1} \log(1+e^{a_rn}) &< \log 2 .
\end{alignat*}
Noting that $a_i > -\varepsilon$, and approximating $\log(1+x) \sim x$, 
this condition becomes
\begin{equation*}
\begin{aligned} & \\ &\Leftrightarrow \end{aligned}
\qquad
\begin{gathered}
(D_1+ \cdots +D_r) n^{r-1} e^{-\varepsilon n} < \log 2 \\
n\varepsilon > (r-1) \log n + 
\log \left (\frac{D_1 + \cdots + D_r}{\log 2} \right).
\end{gathered}
\end{equation*}
If we take $D_1 = \cdots =  D_r = D$ then this simplifies to
$$n\varepsilon > (r-1) \log n + \log \left (\frac{rD}{\log 2} \right).$$
In contrast, Theorem~\ref{thm:hypersurface} tells us that it is
sufficient to take $n$ so that
$$n \varepsilon > (r-1) \log n + \log ((r+3)D) + (r+1)\log 2.$$
Again, this shows that the bounds in Theorem~\ref{thm:hypersurface}
are asymptotically correct, at least for any fixed $r$.  We suspect,
however, that the correct general answer does not have this last term,
or any term which is linear in $r$.


\subsection{Other cyclic resultants}

Instead of the family $\cyclres{f}$, one may wish to consider a more
general family of cyclic resultants.  Let $n_1, \ldots, n_r$ be positive
integers, and consider

$$\gencyclres{f}(\boldz) = 
\prod_{k_1=0}^{n_1-1} \cdots \prod_{k_r=0}^{n_r-1}
f(e^{2 \pi i\ k_1/n_1}\ z_1, \ldots , e^{2 \pi i\ k_r/n_r}\ z_r).
$$

Unfortunately it is not true that the family $\SlopAm_{\gencyclres{f}}$
converges uniformly to $\Am_f$, as $n_1, \ldots, n_r \rightarrow
\infty$.  Trouble occurs if some of the $n_i$ are significantly larger
than others.  For example, consider the amoeba of 
$f(z_1, z_2) = (1+z_1)(1+z_2)$
at a point $(a_1, a_2) \in \RR^2$, with $a_1<0,\ a_2<0$.
Then 
\begin{align*}
\gencyclres{f}(z_1,z_2) &= (1+z_1^{n_1})^{n_2}(1+z_2^{n_2})^{n_1} \\
& = 1 + n_2z_1^{n_1} + n_1 z_2^{n_2} + \cdots
\end{align*}
If $n_2 \sim e^{-a_1 n_1}$, then the first two terms above will have
the same order of magnitude.  Thus $\gencyclres{f}\{(a_1,a_2)\}$ will
not be superlopsided, even if $n_1$ (and hence $n_2$) are large.
It will not be even be lopsided.

However, if we restrict ourselves to the situation in which each
$n_i$ is bounded by some polynomial in each of the other $n_j$, then
an analogous statement to Theorem~\ref{thm:hypersurface} will be
true.  For example, we could let $n_i$ be any polynomial function
of a single parameter $n$.  We will not compute explicit bounds
for approximating the amoeba to within $\varepsilon$, in this more 
general situation; however, the answer will depend upon
these polynomials.  It will certainly still be true that
$\SlopAm_{\gencyclres{f}}$ converges uniformly to $\Am_f$, as this
argument really only depends upon the fact that degrees of
$\gencyclres{f}$ are growing only polynomially, while the terms
are becoming suitably sparse.


\section{Approximating a hypersurface amoeba by linear \\inequalities}
\label{sec:approx}


\subsection{Locating the dominant term}

Theorem~\ref{thm:hypersurface} tells us that for $n$ sufficiently
large, one term of $\cyclres{f}$ dominates, but does not specify which one.
The answer will depend upon in which component of the amoeba 
complement our point $\bolda$ lies.  Since $\cyclres{f} \{\bolda\}$ varies
continuously with $\bolda$, it depends only on the component of the
amoeba complement.

Now the number of components is relatively small compared to the
number of terms of $\cyclres{f}$.  There is a natural 
injective map 
$$\ind :\text{components of $\Amcomp_f$} \hookrightarrow 
\Delta(f) \cap \ZZ^r$$
(c.f. \cite{FPT}).  This is called the
{\bf index} of the component---a complete definition is given below.  
So only a few of the 
terms of $\cyclres{f}$ can possibly be dominant terms.  
Fortunately it is relatively
simple to determine which these are.  The Newton polytope
of $\cyclres{f}$ is $n^r \Delta(f)$, and candidates for dominant 
term are in fact what one would expect them to be: namely, they are the 
images of the integral points of $\Delta(f)$ under this scaling.

\begin{proposition}
Let $\bolda \in \Amcomp_f$, and and let 
$\ind(\bolda) = \veck = (k_1, \ldots, k_r)$
be the corresponding point in $\Delta(f)$.  If 
$\cyclres{f}\{\bolda\}$ is lopsided, then the term of $\cyclres{f}(\boldz)$
which dominates has exponent vector $n^r \veck$, (i.e. it is the
$z_1^{n^r k_1}\cdots z_r^{n^r k_r}$-term).
\end{proposition}

In order to make complete sense of the statement we need to know a definition
of the index $\veck$.  
There are a number of equivalent definitions, but
the simplest for our purposes is the following.

Let $\boldzeta \in \Log^{-1}(\bolda)$.  For each $i \in \{1, \ldots, r\}$,
consider the polynomial,
$$\projf{i}(z) = 
f(\zeta_1, \ldots, \zeta_{i-1}, z, 
\zeta_{i+1}, \ldots, \zeta_r).$$
If $f$ is a polynomial, then
$k_i$ is the number of roots (with multiplicity) of $\projf{i}$
inside the open disc $\{|z| < e^{a_i}\}$.  Since $\projf{i}$ never has
a root on the circle $\{|z| = e^{a_i}\}$ for any $\boldzeta$, this number
is independent of $\boldzeta$.
If $f$ is a Laurent polynomial, then 
\begin{align*}
k_i  &= \#\text{roots of $\projf{i}(z)$ inside $\{0 < |z| < e^{a_i}\}$}
+\mindeg \projf{i}(z) \\
&=\#\text{roots}-\#\text{poles of $\projf{i}(z)$ inside $\{|z| < e^{a_i}\}$}.
\end{align*}

\begin{proof}
It suffices to prove this if $f$ is a polynomial.  If $f$ is a Laurent
polynomial, we multiply by some monomial to make it a polynomial,
and note that this shifts both $\veck$ and the dominant term 
of $\cyclres{f}$ appropriately.

Now 
\begin{align*}
& \projcyclres{f}{i}(z)  \\
&=
\prod_{k_1=0}^{n-1} \cdots \prod_{k_r=0}^{n-1}
f(e^{\frac{2 \pi i k_1}{n}}\ \zeta_1, \ldots, 
e^{\frac{2 \pi i k_{i-1}}{n}}\ \zeta_{i-1},
\ e^{\frac{2 \pi i k_i}{n}}\ z \ ,
e^{\frac{2 \pi i\ k_{i+1}}{n}}\ \zeta_{i+1}, \ldots, 
e^{\frac{2 \pi i k_r}{n}}\ \zeta_r)
\end{align*}
is a product of $n^r$ terms, each of
which is of the form $\projfprime{i}(\phi z)$, for some 
$\boldzeta' \in \Log^{-1}(\bolda)$, $|\phi|=1$.  
Thus the number of roots of $\projcyclres{f}{i}(z)$ inside the
disc $\{|z| < e^{-a_i}\}$ is $n^r k_i$.  As we noted in
Remark~\ref{rmk:domterm},
the number of roots in this disc is actually the exponent of the 
dominant term.

Since the dominant term of $\cyclres{f}$ must project to the
dominant term of each $\projcyclres{f}{i}$, the result follows.
\end{proof}


\subsection{Polyhedral approximations for the components of $\Amcomp_f$}

Fix $n$, and let 
$$\varepsilon = \frac{ (r^2-1) \log n + \log(16/3 cd)}{n},$$
where $c$ and $d$ are defined as in Theorem~\ref{thm:hypersurface}.

Suppose we wish to approximate the component of $\Amcomp_f$ 
corresponding to $\veck \in \Delta(f)$.
From Theorem~\ref{thm:hypersurface}, we know that 
$\SlopAm_{\cyclres{f}} \supset \Am_f$ and approximates $\Am_f$ to 
within $\varepsilon$.  Thus each component of $\Amcomp_f$ is approximated
by a component of $\SlopAmcomp_{\cyclres{f}}$.

Write 
$$\cyclres{f}(\boldz) = M_\veck(z) + \sum_{\vecj \neq n^r \veck}
m_\vecj(\boldz),$$
where 
$$M_\veck(z) = B_\veck z_1^{n^r k_1} \cdots z_r^{n^r k_r}$$
is the
candidate for the dominant term in this component, and 
$m_\vecj(\boldz) = b_\vecj z_1^{j_1} \cdots z_r^{j_r}$.
are the other monomials.

The corresponding component of $\SlopAmcomp_{\cyclres{f}}$ is
the set
$$\Log(\{ \boldz\ \big|\ |M_\veck(\boldz)| > D |m_\vecj(\boldz)|, 
\ \forall\,\vecj\}),$$
where $D+1$ is the number terms in $\cyclres{f}$.
Equivalently, this is the set of $\boldx \in \RR^r$ such that
\begin{equation}
\label{eqn:linearinequalities}
\log |B_\veck| + n^r k_1 x_1 + \cdots + n^r k_r x_r > 
\log D + \log |b_\vecj| + j_1 x_1 + \cdots + j_r x_r
\end{equation}
for all $\vecj$.  This
is a system of linear inequalities in the variable $\boldx$, so 
the solutions to these equations are a convex polyhedron which
approximates the component of the amoeba to within $\varepsilon$.
If there is no component of the $\Amcomp_f$ corresponding to $\veck$
then this system of equations will have no solutions.
Conversely, if this system of inequalities has no solutions,
then this component of the amoeba (if it exists) is not large
enough to contain a ball of radius $\varepsilon$.

Thus we can realise any component of the $\Amcomp_f$ as an increasing
union of convex polyhedra.  This gives an independent proof of the
basic fact \cite{FPT} that the the components of the $\Amcomp_f$ are 
convex. We must admit, however, that there are simpler proofs of
this fact.

Note that in Theorem~\ref{thm:hypersurface}, we actually show 
that $\cyclres{f}$ is $(d n^{r^2 -r})$-superlopsided.
Thus we can in fact take $D = d n^{r^2 -r}$ in 
\eqref{eqn:linearinequalities} and the set of solutions to this
system of inequalities will still
approximate the component of $\Amcomp_f$ to within $\varepsilon$.


\subsection{Approximating the spine}

One of the primary tools for studying amoebas has been the Ronkin
function $N_f$, defined in \cite{Ro}.  For $f \in \basering$, $N_f$
is defined to
be the pushforward of $\log |f|$ under the map $\Log$:
$$N_f(\boldx) :=
\frac{1}{(2\pi i)^r} \int_{\Log^{-1}(\boldx)}
\frac{\log |f(z_1, \ldots, z_r)|\ dz_1 \cdots dz_r}
{z_1 \cdots z_r} .
$$
Ronkin shows in \cite{Ro} that $N_f$ is a convex function, and it is 
affine-linear precisely on the components of $\Amcomp_f$.  When
restricted to a single component of $E$ of $\Amcomp_f$, 
$\nabla N_f = \ind(E)$.

Passare and Rullg{\aa}rd use this function to define the spine of
the amoeba \cite{PR}, as follows. For each component $C$ of $\Amcomp_f$,
extend the locally affine-linear function of $N_f|_E$ to an
affine-linear function $N_E$ on all of $\RR^r$.  Let 
$$N_f^\infty(\boldx) = \max_E \{N_E(\boldx)\} .$$
This is a convex piecewise linear function on $\RR^r$, superscribing
$N_f$. The {\bf spine} of the amoeba $\Am_f$ is defined to be the set
of points where $N_f^\infty$ is not differentiable, and is denoted
$\Spine_f$.

The spine of the amoeba $\Spine_f$ is a strong deformation retract of 
$\Am_f$ \cite{PR, Ru}.
Also note that $\Spine_f$ is actually a tropical hypersurface, as 
defined in the introduction: i.e. it is the singular locus of the maximum
of a finite set of linear functions, where the gradient of each
linear function is a lattice vector.

Now observe that 
$$ \frac{1}{n^r} \log |\cyclres{f}(\boldz)|
=\frac{1}{n^r} \sum_{k_1=1}^n \cdots \sum_{k_r=1}^n 
\log |f(e^{2 \pi i\,k_1/n}\,z_1, \ldots , e^{2 \pi i\,k_r/n}\ z_r)|$$
can be thought of as a Riemann sum for $N_f$.  In particular,
we might expect $\frac{1}{n^r} \log |\cyclres{f}(\boldz)|$ to
converge pointwise to $N_f(\Log(\boldz))$.  This will certainly be
true provided that $\log |\cyclres{f}(\boldz)|$, is bounded on
$\Log^{-1}(\boldx)$, which is the case when 
$\boldx \in \Amcomp_f$.

Suppose $\boldx = \Log(\boldz)$ is in the component of $\Amcomp_f$
of index $\veck \in \Delta(f)$. 
Assume that $\boldx$ has distance
at least $\delta$ from the amoeba, where $\delta>0$ is fixed. 
For any $\varepsilon>0$ we can find $n$ sufficiently large so that
$$\cyclres{f}(\boldz) = M_\veck(\boldz)
+ \sum_\vecj m_\vecj(\boldz),$$
where each $m_\vecj$ is relatively small, i.e.
 $$\sum |m_\vecj(\boldz)| < \varepsilon |M_\veck(\boldz)|$$
(see Corollary~\ref{cor:verylopsided}). Thus we have
$$ 
\log|M_\veck(\boldz)| + \log(1-\varepsilon)
\leq \log|\cyclres{f}(\boldz)| \leq
\log|M_\veck(\boldz)| + \log(1+\varepsilon)
$$
Thus we see that as $n\rightarrow \infty$, the values of 
$\frac{1}{n^r} \log |\cyclres{f}(\boldz)|$,
$$ \frac{1}{n^r}\log|M_\veck (\boldz)|
=\frac{1}{n^r} \log|B_\veck| + k_1 x_1 +  \cdots + k_r x_r
$$
and $N_f(\boldx) = N_E(\boldx) = c_\veck + k_1 x_1 + \cdots + k_r x_r$
all converge upon $N_f(\boldx)$.

We can use this fact obtain good approximations for the spine of
the amoeba.  For each $n$ we consider the function 
$M^\infty:\RR^r \to \RR$ given by
\begin{equation}
\label{eqn:approxspine}
M^\infty(\boldx) := \max \log |M_\veck(\boldz)|,
\end{equation}
where the maximum is taken over all components of $\Amcomp_f$.
This is a piecewise linear function.  We define the approximate spine 
of the amoeba $\LopSpine_{f,n}$ to be the set of points where 
$M^\infty(\boldx)$ is not smooth.  Equivalently $\LopSpine_{f,n}$
is the set of points where the maximum in 
Equation \eqref{eqn:approxspine}
is attained by two distinct values of $\veck$.

\begin{proposition} \ 
\begin{enumerate}
\item $\LopSpine_{f,n} \subset \LopAm_{\cyclres{f}}$.
\item $\lim_{n \to \infty} \LopSpine_{f,n} = \Spine_f$.
\end{enumerate}
\end{proposition}

\begin{proof}
Statement 1 is true because on the component of $\LopAmcomp_{\cyclres{f}}$
of index $\veck$, $|M_\veck(\boldz)| > |M_\vecl(\boldz)|$ for any
other $\vecl \in \Delta(f)$.  Thus the maximum value in Equation
\eqref{eqn:approxspine} cannot be attained by two distinct $\veck$ 
if $\boldx \in \LopAmcomp_{\cyclres{f}}$.

Statement 2 follows from the fact that 
$\frac{1}{n^r} \log|M_\veck(\boldz)| - N_f(\boldx)$ is a constant function
and is less than
$\varepsilon$ for $n$ large.  Let $E_1$ and $E_2$ be components of
$\Amcomp_f$ of index $\veck_1$ and $\veck_2$ respectively.
Consider the hyperplane $H \subset \RR^r$ 
where $\log|M_{\veck_1}(\boldz)|$ and $\log|M_{\veck_2}(\boldz)|$ 
coincide, and the hyperplane $H'$ where $N_{E_1}(\boldx)$ and 
$N_{E_2}(\boldx)$ coincide.  The two hyperplanes $H$ and $H'$
are parallel, and their distance apart is most $\varepsilon K$, 
where $K$ is some constant 
depending only on $\veck_1$ and $\veck_2$.  As there are only
finitely many $\veck \in \Delta(f) \cap \ZZ^r$, these distances can
be made uniformly small.
\end{proof}

For practical reasons, we may wish use an alternate definition
of $M^\infty$, in which one takes the maximum 
in Equation \eqref{eqn:approxspine} over
only those components which appear in $\SlopAmcomp_{\cyclres{f}}$.  
If we do, Statement 1 may be false for small $n$, and we must
settle for saying that $\LopSpine_{f,n} \subset \SlopAm_{\cyclres{f}}$.

One might hope to be able to simplify this construction by taking
the maximum in Equation \eqref{eqn:approxspine} over all
$\veck \in \Delta(f) \cap \ZZ$, rather than just those which
actually correspond to components.  It appears, however, that
this does not give the same answer.
With this alternate definition of $M^\infty$, the `approximate
spine' will have false chambers for all $n$: i.e.  the complement of this 
approximate spine will have components which do not correspond to 
components of $\Amcomp_f$.  We might still hope that these false chambers
would shrink to zero volume as $n$ gets large.
Unfortunately, experimental evidence suggests that the limit 
of these false chambers,
as $n \rightarrow \infty$, can sometimes contain a ball of positive 
radius, and so this method does not produce a good approximation of 
the spine.


\section{More general amoebas}

\subsection{Amoebas of higher codimension varieties in $(\CC^*)^r$}
The higher codimension statement (Theorem~\ref{thm:generalamoeba})
follows fairly quickly from the
hypersurface statement.  Let $V \subset \CC^r$ be a variety, which
is the zero locus of an ideal 
$I = \langle f_1, \ldots, f_k \rangle \subset \basering$.

\begin{proposition}
\label{prop:reduction}
For every $\bolda \in \RR^r$, there exists $f_\bolda \in I$, such that
$$Z_{f_\bolda} \cap \Log^{-1}(\bolda) = V \cap \Log^{-1}(\bolda).$$
\end{proposition}

\begin{proof}
For any Laurent polynomial
$$g(\boldz) = \sum_\vecj b_\vecj z_1^{j_1} \cdots, z_r^{j_r} 
\in \basering,$$ 
let $\bar{g}$ denote its complex conjugate
$$\bar{g}(\boldz) = \sum_\vecj \bar{b}_\vecj z_1^{j_1} \cdots, z_r^{j_r}.$$
We define $f_\bolda$ to be
$$f_\bolda(\boldz) 
:= \sum_{i=1}^k f_i(z_1, \ldots, z_r)
\bar{f}_i(e^{2a_1}z_1^{-1}, \ldots, e^{2a_r}z_r^{-1}).$$
Clearly $f_\bolda$ is a Laurent polynomial and is in $I$. Moreover, 
if we restrict 
$\boldz$ to $\Log^{-1}(\bolda)$,
then $z_i \bar{z}_i = e^{2a_i}$, so
\begin{align*}
f_\bolda(\boldz) &=
 \sum_{i=1}^k f_i(z_1, \ldots, z_r)
\bar{f}_i(\bar{z}_1, \ldots, \bar{z}_r) \\
&= \sum_{i=1}^k f_i(z_1, \ldots, z_r)
\overline{f_i(z_1, \ldots, z_r)} \\
&= \sum_{i=1}^k |f_i(\boldz)|^2.
\end{align*}
Thus $f_\bolda(\boldz) = 0$ if and only if $f_i(\boldz) = 0$ for all $i$.
\end{proof}

This result is also true for ideals in $\CC[z_1, \ldots, z_r]$: one can find
a suitable monomial $m(\boldz)$ such that 
$$m(z_1, \ldots, z_r) 
\bar{f}_i(e^{2a_1}z_1^{-1}, \ldots, e^{2a_r}z_r^{-1}).$$
is a polynomial for all $i$, and a similar argument holds if
$$f_\bolda(\boldz) 
= \sum_{i=1}^k f_i(z_1, \ldots, z_r)
\bigg(m(z_1, \ldots, z_r) 
\bar{f}_i(e^{2a_1}z_1^{-1}, \ldots, e^{2a_r}z_r^{-1})\bigg).$$

As an immediate consequence of Proposition~\ref{prop:reduction} 
we have the following.
\begin{corollary}
\label{cor:intersection}
For any ideal $I \subset \basering$,
$$\Am_I = \bigcap_{f \in I} \Am_f.$$
\end{corollary}

It is also now a simple task to prove our second main result.
\begin{theorem}
\label{thm:generalamoeba}
Let $I \subset \basering$ be an ideal.  A point $\bolda \in \RR^r$ is
in the amoeba $\Am_I$ if and only if $g\{\bolda\}$ is not (super)lopsided
for every $g \in I$.
\end{theorem}

\begin{proof}
If $\bolda \in \Am_I$, then $f\{\bolda\}$ cannot be lopsided for any
$f \in I$, since $\bolda \in \Am_f$ for every $f \in I$.  On
the other hand, suppose $\bolda \notin \Am_I$.  Then by 
Proposition~\ref{prop:reduction} if we take $g = f_\bolda \in I$,
then $\bolda \notin \Am_{g}$.  By Theorem~\ref{thm:hypersurface},
if $n$ is sufficiently large then
$\cyclres{g}\{\bolda\}$ is (super)lopsided,
and $\cyclres{g} \in I$.
\end{proof}

\begin{remark}
\label{rmk:summary}
\rm
In summary, we have three coincident sets, any of which can be used
to define the amoeba $\Am_V$ of a variety $V=V(I)$.
\begin{enumerate}
\item $\Am_V = \Log(V)$.
\item $\Am_V = \bigcap_{f \in I} \Am_f$.
\item $\Am_V = \{ \bolda \in \RR^r\ \big |\ f\{\bolda\}$ is not 
lopsided, for all $f \in I \}$.
\end{enumerate}
In Section~\ref{sec:tropical} we shall see that this is precisely
analogous to a theorem for tropical algebraic varieties.
\end{remark}

If a point $\bolda$ is in $\Amcomp_I$, the proof of 
Theorem~\ref{thm:generalamoeba} also tells us where to look for a witness to 
this fact: namely, we should look at $\widecyclres{f_\bolda}\{\bolda\}$
for all $n$.  For some sufficiently large $n$, this list will be
lopsided.

One unfortunate misfeature of this proof is that it requires us
to use a different $g$ for every point $\bolda \in \Amcomp_I$.  Thus
this statement is purely local.  It does not give any clues as
to how to produce a global uniform approximation to $\Am_I$.  However,
in general we cannot expect there to be any finite set of elements
$g_i \in I$ such that if $\bolda \notin \Am_I$, then some $\widecyclres{g_i}$
is lopsided for $n$ sufficiently large.  If it were so, this would
imply that $\Am_I$ is always an intersection of finitely many
hypersurface amoebas, and this is certainly not true for dimensional
reasons if $\dim V < r/2$.

The best one could for is that we could find some
family of sets $S_n = \{f_{n1}, \ldots, f_{nk_n}\} \subset I$ 
such that 
$$\bigcap_{i=1}^{k_n} \LopAm_{f_{ni}}$$
approximates $\Am_I$ uniformly to within $\varepsilon$ for $n$ 
sufficiently large.  We present it as an open problem to find such 
a family explicitly.


\subsection{Linear transformations and elimination theory}

If we endow $(\CC^*)^r$ with the symplectic form 
$\sum_{i=1}^r \frac{dz_i \wedge \bar{z}_i}{|z_i|^2}$, then the
action of $T=(S^1)^r$ is Hamiltonian, and its moment map is
$\Log$.  Suppose $T'$ is another torus with a 
homomorphism
$\chi: T' \to T$.  Then the induced action of $T'$ on $(\CC^*)^r$
is also Hamiltonian, and we denote its moment map by $\Log'$.  In
this section we describe how to compute the image of 
$V$ under the map $\Log'$.

The image $\Log'(V)$ is always a linear transformation of $\Am_V$.
If $\hat \chi: T^* \to (T')^*$ denotes the adjoint map to $\chi$,
then $\Log' = \hat \chi \circ \Log$
(see e.g. \cite[Prop. 3.2.8]{A}).
More concretely if we 
identify $T'$ with $(S^1)^{r'}$, we can write
$$
\chi(\mu_1, \ldots, \mu_{r'}) =
\bigg( \prod_{i=1}^{r'} \mu_j^{A_{i1}}, \ldots, 
\prod_{i=1}^{r'} \mu_j^{A_{ir}} \bigg),
$$
where $A_{ij}$ are the integer entries of a matrix $A$
---the matrix representation of $\hat \chi$---and
\begin{align}
\Log'(\boldz) &= A\,\Log(\boldz) \label{eqn:lintransform}\\
&= 
\bigg( \sum_{j=1}^r A_{1j} \log|z_j|, \ldots, 
\sum_{j=1}^r A_{r'j} \log|z_j|, \bigg). \notag
\end{align}

We could also take a matrix $A$ with integer entries
as our starting point, 
and construct $T'$, $\chi$, and the map $\Log'$ so
that \eqref{eqn:lintransform} holds.  

Let $I \subset \basering$ denote the ideal of V.  Let $\tilde V
\subset (\CC^*)^{r+r'}$
denote the variety of the ideal 
$$\tilde I = I+J \subset \bigbasering,$$ 
where $J = \langle w_i - \prod_{j=1}^{r} z_j^{A_{ij}} \rangle$.
Now consider the projection of $\tilde V$ onto the $w$-coordinates
$(\CC^*)^{r'}$.  The image of $\tilde V$ under this 
projection is a variety $V'$.  Standard techniques of elimination theory 
allows us to compute its ideal $I'$
(see e.g. \cite{CLO}).

\begin{proposition}
\label{prop:lintransform}
$\Log'(V) = A(\Am_{V}) = \Am_{V'}$.
\end{proposition}

\begin{proof}
A point in $\tilde V$ is a simply pair $(\boldz,\boldw)$ where $\boldz \in V$
and $w_i = \prod_{j=1}^r z_j^{A_{ij}}$.  Thus we have
$$
\Am_{\tilde V} 
= \{(\boldx, \boldy) \in \RR^{r+r'}
\ |\ \boldy = A \boldx,\ \boldx \in \Am_V\}.
$$
Projecting onto the $w$-coordinates, we obtain
\begin{align*}
\Am_{V'} &= 
\{\boldy \in \RR^{r'}\ |\ (\boldx, \boldy) \in \Am_{\tilde V}
\hbox{ for some $\boldx$}\}
\\
&= \{A\boldx \ |\ \boldx \in \Am_V\}  \\
&= A\,\Log(V) \\
&= \Log'(V).
\end{align*}

\end{proof}

It is interesting to note that this construction is closely related
to the cyclic resultants used in the proof of 
Theorem~\ref{thm:hypersurface}.  Suppose $V = Z_f$ is a hypersurface, and 
$\chi:T'=T \to T$ is the map $\chi(t) = t^n$. In this case, $A = nI$ 
is a multiple of the identity matrix, and the variety $V'$
is the zero locus of the function $\cyclres{f}$.  Intuitively, we
should of think that the linear transformation is zooming in on
the amoeba $\Am$; as we zoom in, Theorem~\ref{thm:hypersurface} 
tells us that we see more and more detail in the 
approximations $\LopAm$ and $\SlopAm$.


\subsection{Compactified amoebas}

The most natural generalisation
of amoebas in the compact setting, is to subvarieties of
projective toric varieties.  Each projective toric variety is
a compactification of $(\CC^*)^r$, with an $(S^1)^r$ action which 
extends the $(S^1)^r$ action on $(\CC^*)^r$.  It also carries a 
natural symplectic form $\omega$, for which the $(S^1)^r$ action is
Hamiltonian.  We may therefore use the moment map for this
Hamiltonian action to replace
the map $\Log$.

Our goal in this section is to give a concrete description of this 
more general setting, and observe that our results still hold.
This follows fairly easily from the non-compact case.  Our
construction of toric varieties and their moment maps roughly follows 
a combination of \cite{F} and \cite{A}.

Let $\Delta \subset \RR^r$ be a lattice polytope, i.e. the vertices
of $\Delta$ have integral coordinates.  To every such $\Delta$
we can associate the following data:

\begin{enumerate}
\item A set of lattice points $A = \Delta \cap \ZZ^r$.

\item A semigroup ring $\CC[A]$.  If $A = \{\veck_1, \ldots, \veck_r\}$,
this is defined to be the quotient ring 
$$\CC[s^{\veck_1}, \ldots, s^{\veck_d}]/J,$$
where each $s^{\veck_1}$ has degree $1$, and $J$ is generated by all 
(homogeneous) relations of the form
$$
s^{\veck_{i_1}} \cdots s^{\veck_{i_p}} -
s^{\veck_{j_1}} \cdots s^{\veck_{j_p}} = 0
\hbox{ whenever }
\veck_{i_1}+ \cdots +\veck_{i_p} =  
\veck_{j_1}+ \cdots +\veck_{j_p}
$$
Note that $\CC[A]$ carries an action of the complex torus
$\TT = (C^*)^r$, given by 
$$(\lambda_1, \ldots \lambda_r) \cdot s^{\veck} 
= \lambda_1^{k_1} \cdots \lambda_r^{k_r} s^{\veck}.$$

\item A toric variety $X = \Proj(\CC[A])$.

\item A projective embedding 
$\phi : X \hookrightarrow \PP^{d-1} = \Proj(\CC[t_1, \ldots, t_d])$, 
induced by the map on rings $\CC[t_1, \ldots, t_d] \to \CC[A]$ given
by $t_i \mapsto s^{\veck_i}$.

\item A symplectic form $\omega = \phi^*(\omega_{\PP^{d-1}})$, where
$\omega_{\PP^{d-1}}$ is the Fubini-Study symplectic form on $\PP^{d-1}$.

\item A moment map $\mu$ for the $(S^1)^r$ action on $(X, \omega)$.
We can in fact write down the moment map $\mu$ explicitly.
$$\mu(x) = \frac{1}{\sum_{i=1}^d |s^{\veck_i}(x)|^2}
\sum_{i=1}^d |s^{\veck_i}(x)|^2 \veck_i
$$
To evaluate the right hand side, we must choose a lifting of $x$ to
$\tilde{X} = \Spec(\CC[A])$.  However, since this expression is homogeneous of 
degree $0$ in the $s^{\veck_i}$, 
it is in fact well defined.
\end{enumerate}

It is well known that $\mu(X) = \Delta$, and that if $Y$ is any
other projective toric variety with $\mu_Y(Y) = \Delta$, then 
$Y \cong X$ as toric varieties.

Let $I \subset \CC[A]$ be a homogeneous ideal, and 
$V = \Proj(\CC[A]/I)$ its variety
inside $X$.  
\begin{definition}[Gel'fand-Kapranov-Zelevinsky \cite{GKZ}]
\rm
The {\bf compactified amoeba} of $V$ is $\mu(V) \subset \Delta$.  
We denote the compactified amoeba of $V$ by either
$\overline{\Am}_V$ or $\overline{\Am}_I$ (or by $\overline{\Am}_f$ if
$I = \langle f \rangle$ is principal).
\end{definition}

Let $f \in \CC[A]$ be a homogeneous polynomial of degree $w$.  We
can again decompose $f$ as a sum of monomials, i.e. write
$f = \sum_{i=1}^l m_i$, where each $m_i$ is a $\TT$-weight vector
in $\CC[A]$.  Each of these $m_i$ is a well defined function on
$\tilde{X}$.  Let $\alta \in \Delta$.  We define 
$\overline{f\{\alta\}} := \big\{|m_1(\tilde{\alta})|, \ldots, 
|m_l(\tilde{\alta})|\big\}$ 
where $\tilde{\alta}$ is any preimage of $\alta$ in the composite map
$\tilde{X} \to X \to \Delta$.
Of course, $\overline{f\{\alta\}}$ will depend on the choice of lifting under 
$\tilde{X} \to X$ but only up to rescaling.  Thus the notions of 
$\overline{f\{\alta\}}$ being lopsided or superlopsided are still well defined.
We define $\overline{\LopAm}_f$ and $\overline{\SlopAm}_f$ to
the set of points $\alta \in \Delta$ such that $\overline{f\{\alta\}}$ is 
non-lopsided, and non-superlopsided respectively.

Let $V^\circ$ denote the intersection of $V$ with the open dense
subset of $X$ on which $\TT$ acts freely.  We can identify this
open dense subset of $(\CC^*)^r$, and therefore consider
$\Am_{V^\circ}$.  As both $\Log$ and $\mu|_{(\CC^*)^r}$ are
both submersions with fibres $(S^1)^r$, it follows that
$\Am_{V^\circ}$ is diffeomorphic to $\overline{\Am}_V \cap
\Delta^\circ$, where $\Delta^\circ$ denotes the interior
of $\Delta$.  Let $\psi_\Delta : \Delta^\circ \to \RR^r$ denote
this diffeomorphism.

Moreover any face $\Delta'$ of $\Delta$, corresponds to a toric
subvariety $X' \subset X$.  And $\overline{\Am}_{V} \cap \Delta'
= \overline{\Am}_{V \cap X'}$ \cite{GKZ}.

Thus for every point $\alta \in \Delta$, we can determine
whether $\alta$ is in the compactified amoeba $\overline{\Am}_V$ as
follows.  First we determine the face $\Delta' \subset \Delta$
for which $\alta \in (\Delta')^\circ$.  Then $\psi_{\Delta'}$
identifies $(\Delta')^\circ$ with $\RR^{r'}$ in such a way that
$\overline{\Am}_V \cap (\Delta')^\circ$ is identified
with $\Am_{(V \cap X')^\circ}$. 
We then have $\alta \in \overline{\Am}_V$ if and only if 
$\psi_{\Delta'}(\alta) \in \Am_{(V \cap X')^\circ}$.

\begin{lemma}
\label{lem:uniformtransform}
The map $\psi_\Delta$ is uniformly continuous.
\end{lemma}

\begin{proof}
The projective embedding $\phi$ induces a map
$(\CC^*)^r \hookrightarrow (\CC^*)^{d-1}$, which is define by
monomials.  This induces a linear map from $\RR^r \rightarrow
\Am_{(\CC^*)^r} \subset \RR^{d-1}$.

We also have a map from the moment polytope $\Delta_{d-1}$ of
$(\PP^{d-1})$ to $\Delta$ which is the projection induced by
the inclusion of tori $\TT_{\Delta} \subset \TT_{\Delta_{d-1}}$.

The composite
$$\RR^r \longrightarrow \RR^{d-1} 
\xrightarrow{\psi_{\Delta_{d-1}}}
\Delta_{d-1} \longrightarrow \Delta$$
is $\psi_\Delta$.  Since the first and last maps are uniformly
continuous, it suffices to show that $\psi_{\Delta_{d-1}}$
is uniformly continuous.  

This is fairly straightforward.  Write
$\mu_{\Delta_{d-1}} = (\mu_1, \ldots, \mu_{d-1})$.  For 
$\altz \in (\PP^{d-1})^\circ$ we may write $\altz = (1, z_1, \ldots, z_d)$,
and $\mu_j(\altz) = \frac{|z_j|^2}{1+\sum_{i=1}^{d-1} |z_i|^2}$.
If 
$\big|\log|z_j| - \log |z'_j|\big| < \varepsilon$ then 
$$
\frac{e^{-2\varepsilon}\left(|z'_j|^2\right)}
{e^{2\varepsilon}\left(1+\sum_{i=1}^{d-1} |z'_i|^2\right)} 
< \frac{|z_j|^2}{1+\sum_{i=1}^{d-1} |z_i|^2} < 
\frac{e^{2\varepsilon}\left(|z'_j|^2\right)}
{e^{-2\varepsilon}\left(1+\sum_{i=1}^{d-1} |z'_i|^2\right)} 
$$
i.e.
$e^{-4\varepsilon} \mu_j(\altz') < \mu_j(\altz) < e^{4\varepsilon} \mu_j(\altz')$.
So $|\mu_j(\altz) - \mu_j(\altz)| < (e^{4\varepsilon} -1) \max x_j(\Delta)$.
and the result follows.

\end{proof}

It should therefore come as no surprise that 
Theorems~\ref{thm:hypersurface} and~\ref{thm:generalamoeba} have
analogues in the compact setting.

\begin{corollary}
\label{cor:compacthypersurface}
The families $\overline{\LopAm}_{\cyclres{f}}$ and
$\overline{\SlopAm}_{\cyclres{f}}$ converge uniformly to 
$\overline{\Am}_f$.  Moreover, to approximate the $\overline{\Am}_f$,
to within $\varepsilon$, one can to choose $n>N(\varepsilon)$,
where $N(\varepsilon)$ can be determined solely from $\Delta$ and
$\deg f$.
\end{corollary}

\begin{proof}
One can easily verify that $\overline{g\{\alta\}}$ is (super)lopsided if and
only if $g\{\psi_{\Delta'}(\alta)\}$ is (super)lopsided, where
$\alta \in (\Delta')^\circ$.  Thus, since the lopsidedness
of $\cyclres{f}\{\psi_{\Delta'}(\alta)\}$ can be used to test
membership in $\Am_{f|_{X'}} \subset \RR^{r'}$, 
$\overline{\cyclres{f}\{\alta\}}$ can be used to determine membership in 
$\overline{\Am}_f$.

The fact that the convergence is uniform follows from 
Lemma~\ref{lem:uniformtransform}.  Moreover, as the data $\Delta$
and $\deg f$, are enough to determine an outer bound on the Newton 
polytope of $f$, the rate of convergence can be determined from
these data.
\end{proof}

We also immediately have:
\begin{corollary}
\label{cor:compactamoeba}
$\overline{\Am}_I = \bigcap_{f \in I} \overline{\Am}_f$.  In 
particular $$\overline{\Am}_I = \{\alta \in \Delta\ |\ f\{\alta\}
\hbox{ is not lopsided }\forall f \in I\}.$$
\end{corollary}

It should be noted that Corollary~\ref{cor:compactamoeba} holds
for all toric varieties with a moment map, not just the compact
ones.  However the statement of uniform convergence in 
Corollary~\ref{cor:compacthypersurface} does not 
hold in general for non-compact
toric varieties.  For example, if one considers the toric variety 
$\CC^r$, with the standard
moment map $\mu(\boldz) = \frac{1}{2}(|z_1|^2, \ldots, |z_r|^2)$,
the convergence of the family $\overline{\LopAm}_{\cyclres{f}}$ 
will almost never be uniform.  One can even see this in the simple
example $f(\boldz)= (1-z_1) \cdots (1-z_r)$.  The failure is that 
Lemma~\ref{lem:uniformtransform} does not hold: the map
$\log|x| \mapsto \frac{1}{2}|x|^2$ is not uniformly continuous,
and so the uniform convergence does not carry over.

It is most unfortunate that Proposition~\ref{prop:lintransform} does
not easily carry over to the compact case.  The use of elimination theory
appears only to be well suited to the study of $(\CC^*)^r$ with
its particular standard symplectic form.


\section{Tropical varieties}
\label{sec:tropical}


In this section we show that Theorem~\ref{thm:generalamoeba}
is the analytic counterpart to a theorem for tropical varieties.  
We have already seen examples of tropical
hypersurfaces.  Tropical varieties in general, can be thought
of as a generalisation of amoebas, where one replaces
the norm $|\cdot|:\CC \to \RR$ by a valuation 
in some non-Archimedian field.  For this reason, tropical varieties
are also known as non-Archimedian amoebas.

Let $K$ be an algebraically closed field, with valuation $v$.
For our purposes, a valuation on $K$ is a map 
$v: K \to \Rtrop$, which satisfies the following conditions:
\begin{itemize}
\item $v(xy) = v(x) \odot v(y)$
\item $v(x+y) \leq v(x) \oplus v(y)$.
\end{itemize}
This differs from the usual definition of a valuation in 
two purely cosmetic ways.  First, a valuation is traditionally given as
a map to $v: K \to \RR$; we have simply translated into the
operations of $\Rtrop$.  Second, this is $(-1)$ times the
usual notion of a valuation.  Our reasons for making these
cosmetic changes will become abundantly clear by the end of 
this section.

To every $f \in \Kbasering$ we can associate a tropical polynomial
as follows.  If 
$f = \sum_{\veck \in A} b_\veck z_1^{k_1} \cdots z_r^{k_r}$, write
\begin{align*}
f_\tau(\boldx) &= \bigoplus_{\veck \in A} v(b_\veck)\odot \boldx^\veck \\
& = \max_{\veck \in A} \{v(b_\veck) + \boldx \cdot \veck\}
\end{align*}
and call it the {\bf tropicalisation} of $f$.  
We denote the tropical hypersurface
associated to $f_\tau$ by $\TropVar_f$.

If $\bolda \in \Rtrop^r$, we can assign 
a weight to every monomial $m \in \Kbasering$:  define the weight
of $m$ at $\bolda$ to be $$\wt_\bolda(m) := m_\tau(\bolda).$$
If $f(\boldz) = \sum_{i=1}^d m_i(\boldz)$ where $m_i$ are
monomials, let
$$f\{\bolda\}_\tau = \{\wt_\bolda(m_1), \ldots, \wt_\bolda(m_d)\}.$$
Recall that in $\Rtrop$, a list of numbers $\{b_1, \ldots,  b_r\}$
is (tropically) lopsided if the maximum of element of this list does
not occur twice (in which case the maximum element is greater
than the tropical sum of all the other elements).  Thus 
$f\{\bolda\}_\tau$ is lopsided if and only if $\bolda \notin \TropVar_f$.

Let $I \subset \basering$ be an ideal, and
$$V = V(I) = \{\boldz \in (K^*)^r\ |\ f(\boldz)=0\ \forall f \in I\}$$
be its affine variety.  Let
$\val: (K^*)^r \to \Rtrop^r$ be the map
$$\val(\boldz) = (v(z_1), \ldots, v(z_n)).$$

The following theorem, as stated, most closely resembles the formulation
in \cite{SS}, though variants of it have also appeared in 
\cite{EKL, St}. 

\begin{theorem}[Speyer-Sturmfels \cite{SS}]
\label{thm:tropical}
The following subsets of $\Rtrop^r$ coincide:
\begin{enumerate}
\item The closure of the set $\val(V)$.
\item The intersection of all tropical hypersurfaces 
$\bigcap_{f\in I} \TropVar_f$
\item The set of points $a \in \Rtrop^r$, such that $f\{\bolda\}_\tau$
is not lopsided for all $f \in I$.
\end{enumerate}
\end{theorem}
This set is called the {\bf tropical variety} of the ideal $I$.

In fact a stronger result than Theorem~\ref{thm:tropical} 
(as stated here) is shown in
\cite{SS}.  Let $k$ denote the residue
field of $K$.  If $I \subset K[z_1, \ldots, z_r]$,
then one can construct an initial ideal of I, in $k[z_1, \ldots, z_r]$,
corresponding to any weight $\bolda \in \Rtrop$.  One can equivalently
describe the tropical variety of $I$ as the set of points $\bolda \in \Rtrop$
such that the associated initial ideal contains no monomial.  Thus
it is possible to determine membership in a tropical variety using
Gr\"obner basis techniques.

One can easily see that Theorem~\ref{thm:tropical} is precisely analogous to
the summary given in Remark~\ref{rmk:summary}.  The proofs of
these results, however, are extremely different.  An obvious
question, therefore, is whether analogous statements can be made in
other contexts.  

The following is a general context in which one might hope for
such a theorem to be true.
Suppose that $K$ is an algebraically closed field, and let 
$S(\odot, \oplus, \leq)$ be a totally ordered semiring.  
Suppose $\|\cdot\|_K: K^* \to S$ satisfies
the following conditions:
\begin{enumerate}
\item $\|xy\| = \|x\|_K \odot \|y\|_K$ for all $x,y \in K$.
\item for all $a,b \in S$ we have 
$$a \oplus b = \max \big\{\|x+y\|_K\ \big|\ \|x\|_K=a,\ \|y\|_K=b \big\}.$$
\end{enumerate}
In particular condition 2 implies that 
$\|x+y\|_K \leq \|x\|_K \oplus \|y\|_K$ for all $x, y \in K$.  Thus 
$\|\cdot\|_K$ is an $S$-valued norm.

Let $f \in  \Kbasering$, and write $f = \sum_{i=1}^d m_i$
as a sum of monomials.  For any point $\bolda \in S^r$, let 
$\boldzeta$ be such that $\|\boldzeta\|_K = \bolda$. 
We define
$$f\{\bolda\} := \big\{
\|m_1(\boldzeta)\|_K, \ldots, 
\|m_d(\boldzeta)\|_K \big\}.$$
As $\|\cdot\|_K$ is multiplicative, this is independent of the
choice of $\boldzeta$.  Since $S$
is totally ordered, we can define a list to of elements of
$S$ to be {\em lopsided} if and only if one number is greater
than the sum of all the others.

Let $V \subset (K^*)^r$ be a variety defined by an ideal
$I \subset \Kbasering$.  We can consider the following
sets:
\begin{itemize}
\item The closure of 
$\big\{(\|z_1\|_K, \ldots, \|z_r\|_K)\ \big|\ \boldz \in V\big\}$
\item $\{\bolda \in S^r\ |\ f\{\bolda\}$ is not lopsided
$\forall f \in I\}$
\end{itemize}

The question is whether these two sets are equal for a particular
$(K, S, \|\cdot\|_K)$.
In this paper, we primarily discussed the example in which $K = \CC$,
$S = \RR_+$, and $\|\cdot\|_\CC = |\cdot|$, and showed that they are 
equal.  We have 
also just seen that this is true if $K$ is a non-Archimedian field, 
with $\|\cdot\|_K$ as its valuation and $S = \Rtrop$.  

Many (though
not quite all) of the elements
of the proof of Theorem~\ref{thm:hypersurface} are valid in a more
general context.
Suppose that in addition to being a totally ordered
semiring, $S$ is a $\QQ_+$-module (i.e. we can make sense of
such things as $\frac{2}{3}a$ for $a \in S$); for example
$\Rtrop$ is a $\QQ_+$-module with the trivial $\QQ_+$-action.

Define a binary operation $\ominus$ on $S$ by
$$a \ominus b := \min\ \big\{c \in S\ \big|\ c \oplus b \geq a \big\}$$
whenever this set is nonempty. (We need not overly concern ourselves 
with
the fact that a precise minimum may not exist: one can always 
get around this by treating this set as a Dedekind cut.)
Then the triangle inequality $\|x-y\|_K \geq \|x\|_K \ominus \|y\|_K$
is valid (assuming $\|x\|_K \geq \|y\|_K$).  To see this note that
$a \leq a'$ implies $a\ominus b \leq a'\ominus b$; thus
$$\|x\|_K \ominus \|y\|_K 
\leq  \big(\|x-y\|_K \oplus \|y\|_K\big) \ominus \|y\|_K.$$
Clearly $\|x-y\|_K \in \big\{c \in S\ \big|\ c \oplus \|y\|_K 
\leq \|x-y\|_K \oplus \|y\|_K \big\}$
which implies that 
$$\|x-y\|_K \geq \big(\|x-y\|_K \oplus \|y\|_K\big) \ominus \|y\|_K.$$

A closer examination of the proofs
of Lemma~\ref{lem:goodestimate} and
Calculation~\ref{cal:expderiv} now reveal that they are also valid 
(almost word for word) for a general $(K,S,\|\cdot\|_K)$.
We can also prove Lemma~\ref{lem:averaging} in general, by replacing
the integral over the torus
$$\frac{1}{(2\pi i)^r}
\int_{|z_1|=1} \cdots \int_{|z_r|=1}
\left(
\sum_\vecj 
\frac{m_\vecj (e^{a_1}z_1, \ldots, e^{a_r}z_r)}
{z_1^{l_1} \cdots z_r^{l_r}}
\right)
\frac{d z_1}{z_1} \cdots \frac {d z_1}{z_1}
$$
by a discrete average over a finite subgroup of the torus
$$
\frac{1}{N^r}
\sum_{z_1:z_1^N=1} \cdots \sum_{z_r:z_r^N=1}
\left(
\sum_\vecj 
\frac{m_\vecj (e^{a_1}z_1, \ldots, e^{a_r}z_r)}
{z_1^{l_1} \cdots z_r^{l_r}}
\right).
$$
If $N$ is suitably large this discrete average will have the same
effect as the integral (i.e. picking out a single term from the
polynomial).  In fact we can follow the proof Statement 1 of
Theorem~\ref{thm:hypersurface}, up to and including the inequality 
\eqref{eqn:b}.  All that
remains to show that the right hand side of \eqref{eqn:b} becomes
sufficiently small as $n$ gets large.  Unfortunately, in general,
this is not always true for all $\gamma < 1$
(e.g. if $S = \RR[[x]]$ with $f < g$ if the leading coefficient
of $f-g$ is positive, this fails for $\gamma = 1-x$).
It is true, however, for $S=\Rtrop$, so our proof of 
Theorem~\ref{thm:hypersurface} actually works in the tropical setting
as well. This suggests that it may be possible to prove a
more general statement, with perhaps some additional restrictions
on $S$.  
 
On the other hand, our proof of Corollary~\ref{cor:intersection}, 
which is the
key to Theorem~\ref{thm:generalamoeba}, is very highly dependent
on the fact that our underlying field is the complex numbers.
A very different proof is required in the tropical setting, and
it would be interesting to try and unite these two.


\appendix

\section{Appendix: Details of calculations}
\label{app:calcdetails}


In this appendix we give the mundane, technical calculations
which are used in the proofs of Lemma~\ref{lem:onevariable}
and Theorem~\ref{thm:hypersurface}.

\begin{calculation}
\label{cal:1}
If $n \log \gamma^{-1} \geq (D_0 + D_1)\log n + \log (8/3 c_0 c_1)$, then
\begin{enumerate}
\item $\mdeg{1} ((1 + \gamma^n)^{\mdeg{0}} - 1) < 1/2$,
\item $(1+ \gamma^n)^{\mdeg{0}} < 3/2$.
\end{enumerate}
\end{calculation}

\begin{proof}

\begin{alignat}{2}
&& n \log \gamma^{-1} &\geq {D_0 + D_1} \log n + \log (8/3 c_0 c_1) \notag\\
&\Leftrightarrow\qquad & \gamma^{-n} &\geq 8/3 \mdeg{0} \mdeg{1} \notag\\
&\Leftrightarrow\qquad & \mdeg{0} \gamma^n &\leq \frac{3}{8 \mdeg{1}} .
   \label{eqn:app1}
\end{alignat}

Also
\begin{align}
\frac{3}{8 \mdeg{1}} &\leq 
\left(\frac{1}{2 \mdeg{1}}\right) 
- \frac{1}{2}\left(\frac{1}{2 \mdeg{1}} \right )^2 \notag\\
& < \log \left (1+ \frac{1}{2 \mdeg{1}} \right ) \label{eqn:app2} \\
& < \log (3/2) . \label{eqn:app3}
\end{align}

Using \eqref{eqn:app1}, \eqref{eqn:app2} and the fact that
$\log (1+\gamma^n) < \gamma^n$, we find
\begin{alignat*}{2}
&& \mdeg{0} \log(1+\gamma^n) &
< \log \left (1 + \frac{1}{2 \mdeg{1}} \right ) \\
&\Rightarrow\qquad & (1+\gamma^n)^{\mdeg{0}} &< 1+ \frac{1}{2 \mdeg{1}} \\
&\Rightarrow\qquad & \mdeg{1}((1 + \gamma^n)^{\mdeg{0}} - 1) &< 1/2 .\\
\end{alignat*}

On the other hand from \eqref{eqn:app1} and \eqref{eqn:app3} we have
\begin{alignat*}{2}
&& \mdeg{0} \log(1+\gamma^n) &< \log (3/2) \\
&\Rightarrow\qquad & (1+ \gamma^n)^{\mdeg{0}} &< 3/2 .
\end{alignat*}

\end{proof}

\begin{calculation}
\label{cal:lopbounds}
If
$n\log \gamma^{-1} \geq (r-1) \log n + \log \left((r+3)2^{r+1}c\right)$,
then
$$ 
\frac{ e^{(r+2)c n^{r-1} \gamma^n} -1}
{2-e^{c n^{r-1} \gamma^n}}
< 2^{-r-1}.
$$
\end{calculation}

\begin{proof}
\begin{alignat}{2}
&& n\log \gamma^{-1} &\geq (r-1) \log n + \log \big((r+3)2^{r+1}c\big) 
\notag\\
&\Leftrightarrow\qquad & \gamma^{-n} &\geq  (r+3)2^{r+1} c n^{r-1}\notag\\
&\Leftrightarrow\qquad & c n^{r-1}\gamma^n &\leq \frac{2^{-1-r}}{r+3}  .
\label{eqn:app4}
\end{alignat}

Using the power series expansion of $\log(1+x)$, one can
check that
\begin{equation}
\label{eqn:app5}
(r+2)\frac{2^{-1-r}}{r+3} < 
\log \left(1 + \frac{(2+r)2^{-1-r}}{2+r+2^{-1-r}} \right)
\end{equation}
for all $r \in \ZZ_+$.  Thus from \eqref{eqn:app4} and \eqref{eqn:app5}
we have:
\begin{alignat}{2}
&& (r+2) c n^{r-1}\gamma^n &\leq 
\log \left(1 + \frac{(2+r)2^{-1-r}}{2+r+2^{-1-r}} \right) \notag \\
&\Leftrightarrow\qquad & e^{(r+2) c n^{r-1} \gamma^n} &\leq  
1 + \frac{(2+r)2^{-1-r}}{2+r+2^{-1-r}} \notag \\
&\Leftrightarrow\qquad & e^{(r+2) c n^{r-1} \gamma^n} -1 &\leq  
 \frac{(2+r)2^{-1-r}}{2+r+2^{-1-r}} . \label{eqn:app6} 
\end{alignat}

Similarly one can check that
\begin{equation}
\label{eqn:app7}
\frac{2^{-1-r}}{r+3} <
\log \left(1 + \frac{2^{-1-r}}{2+r+2^{-1-r}} \right)
\end{equation}
for all $r \in \ZZ_+$.  So from \eqref{eqn:app4} and \eqref{eqn:app7}
we have:
\begin{alignat}{2}
&& c n^{r-1}\gamma^n &\leq 
\log \left(1 + \frac{2^{-1-r}}{2+r+2^{-1-r}} \right) \notag\\
&\Leftrightarrow\qquad & e^{c n^{r-1} \gamma^n} &\leq  
1 + \frac{2^{-1-r}}{2+r+2^{-1-r}} \notag \\
&\Leftrightarrow\qquad & 2- e^{c n^{r-1} \gamma^n} &\geq  
1 - \frac{2^{-1-r}}{2+r+2^{-1-r}} \notag \\
&& &= \frac{2+r}{2+r+2^{-1-r}} . \label{eqn:app8}
\end{alignat}

Putting together \eqref{eqn:app6} and \eqref{eqn:app8}:
\begin{align*}
\frac{ e^{(r+2) c n^{r-1} \gamma^n} -1}
{2-e^{c n^{r-1} \gamma^n}}
&< \frac{(2^{-1-r})\,(2+r)/(2+r+2^{-1-r})}
{(2+r)/(2+r+2^{-1-r})} \\
&= 2^{-1-r} .
\end{align*}

\end{proof}

\begin{calculation}
\label{cal:expderiv}
For $x>0$, and $s \in \ZZ_+$,
$$
\sum_{w_0 \geq 1} {w_0 + s-1 \choose s-1}
\sum_{w \geq w_0 }
\frac{x^w}{w!}
<  e^{(s+1)x} -1
$$
\end{calculation}

\begin{proof}
\begin{align*}
\sum_{w_0 \geq 1} {w_0 + s-1 \choose s-1}
\sum_{w \geq w_0 }
\frac{x^w}{w!}
&< -1 + \sum_{w_0 \geq 0} {w_0 + s-1 \choose s-1}
\sum_{w \geq w_0 }
\frac{x^w}{w!} \\
&= 
-1 + \sum_{w \geq 0} \sum_{w_0 =0}^w {w_0 + s-1 \choose s-1}
\frac{x^w}{w!} \\
&= 
-1 + \sum_{w \geq 0}  {w + s \choose s}
\frac{x^w}{w!} \\
&= -1 + \sum_{w \geq 0}  \frac{1}{s!} \frac{d^s}{dx^s}
\left(
\frac{x^{w+s}}{w!}\right) \\
&= -1 + \frac{d^s}{dx^s}
\left(
\frac{x^s}{s!}
\sum_{w \geq 0}  
\frac{x^w}{w!}\right) \\
&= -1 + \frac{d^s}{dx^s} \left( \frac{x^s}{s!} e^x \right) \\
&= -1 + \left(\sum_{i=0}^s {s \choose i}\frac{x^i}{i!}\right) e^x \\
&< -1 + \left(\sum_{i \geq 0} \frac{(sx)^i}{i!}\right) e^x \\
&< -1 + e^{(s+1)x} .
\end{align*}

\end{proof}


\end{document}